\documentclass[12pt]{article}

\usepackage[utf8]{inputenc} 
\usepackage[T1]{fontenc}    
\usepackage{hyperref}       
\usepackage{url}            
\usepackage{booktabs}       
\usepackage{amsfonts}       
\usepackage{nicefrac}       
\usepackage{microtype}      
\usepackage{lipsum}		
\usepackage{graphicx}
\usepackage{natbib}
\usepackage{doi}

\usepackage{natbib}
\bibpunct[, ]{(}{)}{,}{a}{}{,}%
\usepackage{multirow}
\usepackage{graphicx}
\usepackage{array}
\usepackage{longtable}
\usepackage{enumitem}
\usepackage{placeins}
\usepackage{tabularx}
\usepackage{csquotes}

\newlist{steps}{enumerate}{1}
\setlist[steps]{label=\textit{Step \arabic*:},leftmargin=*}
   
\usepackage{algorithm}
\usepackage{algorithmic}
\usepackage{setspace}
\usepackage{easyReview}
\usepackage{clipboard}
\usepackage{amsmath,amssymb}
\DeclareMathOperator{\EX}{\mathbb{E}}
\DeclareMathOperator{\PX}{\mathbb{P}}

\usepackage{amsthm} 
\newcounter{propcount}
\newtheorem{proposition}[propcount]{Proposition} 

\newcounter{theocount}
\newtheorem{theorem}[theocount]{Theorem} 

\DeclareMathOperator*{\argmin}{arg\,min} 

\usepackage{titlesec}
\titlespacing{\section}{0pt}{8pt}{0pt}
\titlespacing{\subsection}{0pt}{8pt}{0pt}
\titlespacing{\subsubsection}{0pt}{8pt}{0pt}

\usepackage{caption}
\usepackage[printonlyused,nohyperlinks]{acronym}
\acrodef{MDP}{\textit{Markov decision process}}
\acrodefplural{MDP}{\textit{Markov decision processes}}
\acrodef{ADP}{\textit{approximate dynamic programming}}
\acrodef{ALP}{\textit{approximate linear program}}
\acrodef{OR}{\textit{operations research}}
\acrodef{HRSP}{\textit{home care routing and scheduling problem}}
\acrodef{VRP}{\textit{vehicle routing problem}}
\acrodef{PVRP}{\textit{periodic vehicle routing problem}}
\acrodef{LP}{\textit{linear program}}
\acrodef{1D}{\textit{one-dimensional}}
\acrodef{2D}{\textit{two-dimensional}}

\usepackage{geometry}
\title{Dynamic Home Care Routing and Scheduling with Uncertain Number of Visits per Referral}

\author{\large Danial Khorasanian$^{}$\thanks{Corresponding Author, Email: d.khorasanian@gmail.com, Danial.Khorasanian@telfer.uottawa.ca}, Jonathan Patrick$^{}$\thanks{Email: Patrick@telfer.uottawa.ca}, Antoine Saur\'{e}$^{}$\thanks{Email: asaure@uottawa.ca} \\ \\
{\footnotesize Telfer School of Management, University of Ottawa, Ottawa, Ontario K1N 6N5, Canada} \\ \\
\textsc{\footnotesize A Preprint - March 12, 2024}} 

\date{} 

\usepackage{fancyhdr}
\usepackage{lipsum} 
\usepackage{titlesec}

\titleformat*{\section}{\large\bfseries}
\titleformat*{\subsection}{\normalsize\bfseries}
\titleformat*{\subsubsection}{\small\bfseries}

\begin{document}

\maketitle

\doublespacing
\begin{abstract}
Despite the rapid growth of the home care industry, research on the scheduling and routing of home care visits in the presence of uncertainty is still limited. This paper investigates a dynamic version of this problem in which the number of referrals and their required number of visits are uncertain. We develop a \textit{Markov decision process} (MDP) model for the single-nurse problem to minimize the expected weighted sum of the rejection, diversion, overtime, and travel time costs. Since optimally solving the MDP is intractable, we employ an \textit{approximate linear program} (ALP) to obtain a feasible policy. The typical ALP approach can only solve very small-scale instances of the problem. We derive an intuitively explainable closed-form solution for the optimal ALP parameters in a special case of the problem. Inspired by this form, we provide two heuristic reduction techniques for the ALP model in the general problem to solve large-scale instances in an acceptable time. Numerical results show that the ALP policy outperforms a myopic policy that reflects current practice, and is better than a scenario-based policy in most instances considered. \vspace{3pt}


\end{abstract}


\section{Introduction}\label{sec:introduction}

Home care refers to medical, paramedical, and social services provided
for patients at their home \citep{Cappanera2018}. Cost advantages
and an aging population have made home care one of the world's fastest growing
industries \citep{Heching2019}. The \ac{HRSP} is one of the most common optimization problems arising
in this area \citep{Fikar2017}. 

In the home health care setting, each new request represents a referral with one or, more likely, several visits with a specific pattern (e.g., one visit every three days or visits on Mondays and Thursdays). The \ac{HRSP} is a joint decision problem that includes accepting or rejecting new referrals, scheduling
the accepted referrals, assigning them to eligible nurses, and
routing the daily visits for each nurse. Poor solutions to this problem
not only affect the quality of the services provided to home care
patients and the cost incurred by home health care providers but also
result in delays in the discharge of hospital patients whose health
needs could be met by home care services. This, in turn, leads to
additional congestion in hospitals and therefore delays in medical
interventions requiring hospital stays.

This research was inspired by a collaboration of the authors with a not-for-profit home health care organization in Canada. One of their primary challenges in solving the \ac{HRSP} is the uncertainty in the number of new referrals and the number of visits needed for each referral. The uncertainty in the number of visits is the result of a number of factors. For example, the speed of recovery in certain pain management or post-surgery care scenarios is unknown as it depends on the patient's features and the nurses' skills \citep{Ferrari2018}. Alternatively, a patient may need to return to the hospital unexpectedly cutting short his/her home care services, or the type of service for a patient may change abruptly requiring a different nurse with a different skill set. Thus, we investigate a dynamic HRSP with a stochastic number of referrals per day and a stochastic number of visits for each referral. We develop our model for a single nurse, and briefly describe how our work can be extended to provide a reasonable policy in the multi-nurse setting.

Our model captures some unique information about the referrals such as the maximum number of days a new referral can wait for the first visit based on her/his urgency level, and the probability distribution for the number of visits. The nurse has a limited shift length each day. Once a referral is accepted by the provider, all of her/his visits must be served until discharge. If this results in an excessive workload on a given day, two options are available to the provider. First, the working hours can be extended through (limited) overtime at a cost. Alternatively, any visit assigned to a day can be diverted to an on-call nurse at a higher cost. It is necessary to consider diversions in our problem because the number of visits is uncertain, rejection in the middle of care is impossible, and overtime is limited. The cost function in our model is a weighted sum of the rejection, diversion, overtime, and travel time costs.

We formulate the problem as a \ac{MDP} and use an \ac{ALP} to overcome the curse of dimensionality. The \ac{ALP} reduces the computational burden by approximating the value function of the MDP using an affine parametric form and acting greedily with respect to it. In comparison with the simulation-based \ac{ADP} approaches, the ALP approach is not subject to simulation error and has more stable and monotone properties \citep{Adelman2004}. However, even the ALP approach runs into computational challenges for large-scale instances. We present a closed-form solution for the optimal \ac{ALP} parameters of a special case of the problem. We then propose two heuristic reduction techniques for the ALP model of the general problem by observing certain features of the optimal \ac{ALP} parameters for the special case. These techniques enable us to find good policies for large-scale instances of the problem in an acceptable time. We derive two benchmark policies to compare against our ALP policy and obtain an estimated upper bound on the optimality gap of the ALP policy for a specific instance using a complete information relaxation.

\subsection{Our contribution}\label{sec:our-cont}
We develop an ALP approach to solve a single-nurse dynamic HRSP with a stochastic number of referrals per day and a stochastic number of visits per referral. We use a novel strategy to make this approach tractable for our problem using two heuristic reduction techniques based on a closed-form solution for a special case of the problem. Our ALP policy outperforms the benchmark policies in the majority of the scenarios we examined, especially in the most general case of our problem, and offers valuable practical insights that home health care providers can utilize for capacity planning. Our model, developed through conversations with our partner provider, is more comprehensive than those in previous studies of the dynamic HRSP as it addresses the stochasticity of the number of visits per referral and includes additional features, such as the capability to divert visits and use overtime, and the flexibility to set distinct service times, care patterns and wait-time targets for various service types. These additional features increase the realism of our model and provide more flexible solutions that may reduce the total cost.

\subsection{Related literature}\label{subsec:Related-literature}
The literature on the \ac{HRSP} has been growing very fast in recent years. A comprehensive review outlining different variations of the problem can be found in \citep{Fikar2017}. However, the literature regarding stochastic \ac{HRSP}s is still limited.

\cite{Cire2022} proposed a one-step policy improvement \ac{ADP} approach for the dynamic single-period \ac{HRSP} with multiple nurses. The objective is to determine a policy that maximizes the expected discounted reward over an infinite horizon, where the reward function is based on the revenue from accepting referrals, the cost of rejection, and the cost of nurses’ wages. They demonstrate the superiority of their method over some myopic approaches. Unlike our model, theirs assumes a deterministic number of visits per referral, and these visits occur exclusively on consecutive days. Moreover, it lacks the ability to divert visits, assumes a fixed service time of 45 minutes per visit for all referrals, and restricts the first visit of any accepted referral to be assigned to the current day only.

Among the rare studies in the literature on the dynamic HRSP, only \cite{Cappanera2018} considered the uncertainty in the number of visits in their modeling but with a different problem setting and objective compared to ours and over a short-term finite horizon as opposed to the infinite horizon in our case. They present a non-standard cardinality-constrained robust approach for the multi-nurse \ac{HRSP} with uncertain new arrivals and current patients' demand. They assume that each patient needs a fixed number of visits with certainty and some uncertain visits in the planning horizon. For example, a patient definitely needs three visits in a week, and may need two additional visits in the same week. Each nurse's daily tour includes at most a given number of uncertain visits. The objective is to balance the workload of nurses while considering the maximum number of uncertain requests.

\cite{Bennett2011} investigated a dynamic single-nurse \ac{HRSP}. They developed some rolling horizon myopic heuristics to maximize
the number of visits served in the planning horizon. Decisions concerning new arrivals follow a greedy approach that ignores uncertain future arrivals. Each patient needs service for a given number of weeks with a given number of visits per week (e.g., 8:00 AM on Mondays and 10:00 AM on Wednesdays over four weeks). They mention in their conclusion that imposing a fixed time for all visits of a patient is quite rigid and that more flexibility would be desirable. In a follow-up paper, \cite{Demirbilek2019a} proposed a scenario-based heuristic in order to consider unknown future demand in the problem described by \cite{Bennett2011}. Then, in a second paper, \cite{Demirbilek2019b} extended this method to solve a multi-nurse version of the problem. All three studies assume a fixed appointment time for all visits of each referral, and a fixed service time per visit equal to 30 minutes for all referrals. 

The HRSP is an application of the well-known \ac{VRP}, and our problem can be characterized as a variant of the \ac{PVRP}. In the \ac{PVRP}, customers need one or more visits during a planning horizon of several periods and with a specific pattern (Rothenbacher 2019). Among the applications of the \ac{PVRP} are preventive maintenance, delivery of blood products to hospitals, and waste materials collection \citep{Rodriguez-Martin2019}. The most recent literature review on the \ac{PVRP}, \citep{Campbell2014}, confirms that the literature on the stochastic \ac{PVRP} is very limited. They only mention \cite{Schedl2011350} who studied a PVRP in which the customer's demand at each visit is stochastic. This problem is modeled as a stochastic integer program with recourse and is solved using some metaheuristics. To the best of our knowledge, there are no studies on the PVRP that consider a stochastic number of visits.

To the best of our knowledge, we present the first study on the application of the \ac{ALP} approach to the dynamic \ac{HRSP} or the \ac{PVRP} in the literature. Among the most foundational research on the ALP approach are the papers by \cite{SCHWEITZER1985568,DeFarias2003850}, and \cite{Adelman2004}. This approach has been used successfully to obtain feasible policies for dynamic scheduling problems such as dynamic multi-priority patient scheduling \citep{Patrick2008} and dynamic multi-appointment patient scheduling \citep{Saure2012}. Moreover, this approach or some hybrid versions of it have produced excellent feasible policies for variations of the VRP involving an inventory/routing problem with stochastic demand \citep{Adelman2004} and a single-vehicle dynamic dispatch wave problem with one-dimensional distances \citep{Klapp2016}. In addition to generating a feasible policy, the ALP approach can be used to provide a bound on the optimal value function, and this bound can be applied to obtain optimality gaps for heuristic policies (see, e.g., \citealp{Adelman2004}). \cite{DeFarias2003850} outline the issues with this bound and suggest some ideas to alleviate them. 

Despite the advantages of the ALP approach, it is sometimes complex to solve due to the often large number of constraints and/or variables. In these cases, some techniques must be applied to reduce the size of the model. For example, \cite{Tong2014121} and \cite{Vossen20151352} showed that the original ALP for the network revenue management problem can be represented by a more compact \ac{LP} that can be solved much faster.

The primary distinguishing characteristic of our problem, in contrast to the stochastic TSP or VRP, is that requests involving multiple visits should be served according to a specific pattern across distinct days, rather than a single visit in a day. Among recent applications of ADP methods to a VRP with some similarities to our problem, \cite{Ulmer2019185} combined an offline value function approximation with an online rollout algorithm to solve a dynamic VRP with stochastic requests in a single period. In this study, the decision maker assigns or rejects each new request and directs the vehicle to the next location at each decision epoch.

\subsection{Paper structure}
The rest of the paper is organized as follows. In Section \ref{sec:problem_des}, we describe the problem. In Section \ref{sec:MDP_formulation}, we formulate the \ac{MDP} model. In Section \ref{sec:Approximate-linear-program}, we formulate the \ac{ALP} approach. In Section \ref{sec:A-closed-form-solution}, we determine a closed-form solution for the \ac{ALP} parameters in a special case of the problem. In Section \ref{sec:Acceleration-techniques-for}, we present two heuristic reduction techniques to be used in the ALP approach for the general problem. In Section \ref{sec:Numerical-results}, we evaluate the performance of the ALP policy through numerical experiments. In Section \ref{sec:Conclusions-and-suggestions}, we state our conclusions and suggestions for future research. In Appendix A, we describe our benchmark policies in detail. In Appendix B, we provide a discussion on lower bounds for our MDP. In Section \ref{sec:ALP_tuning} of the electronic companion, we explain the procedure to tune the ALP. In Section \ref{sec:proof-theorem1}, we prove the closed-form solution presented for the special case of the problem. In Section \ref{sec:effects_of_accel}, we analyze the heuristic reduction techniques for the ALP through some numerical experiments. In Section \ref{sec:E_add_num}, we report some additional information about the numerical experiments. Finally, we develop a simple heuristic to solve the multi-nurse version of the problem in Section \ref{sec:E_multi-nurse}.  

\section{Problem description}\label{sec:problem_des}
\begin{figure}[h]
	\begin{center}
		\includegraphics[height=1.7in]{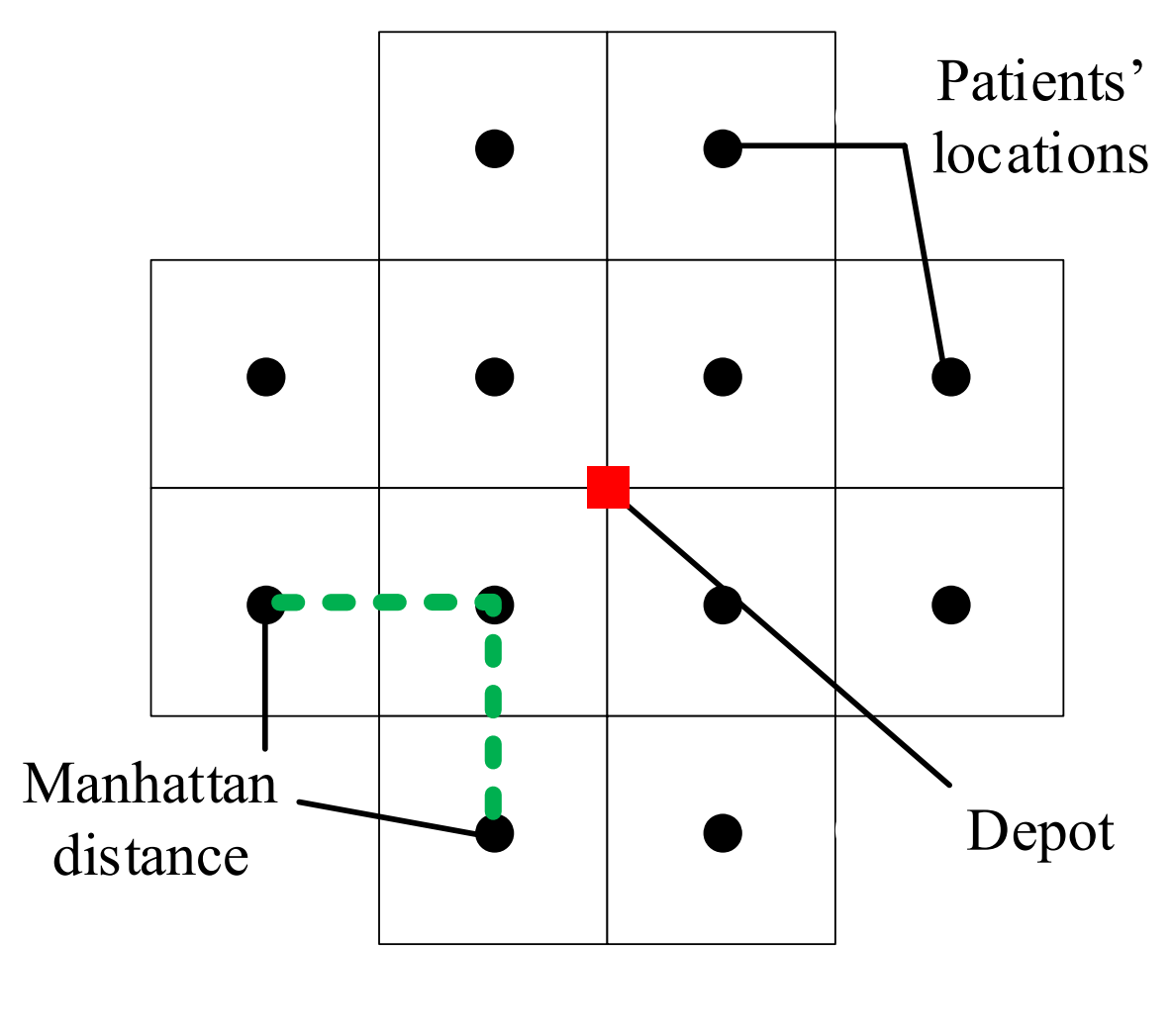}
		\caption{A circular grid area with Manhattan distances around the depot} \label{grid}
	\end{center}
\end{figure}

A geographical area represented by a circular grid-shaped area (Figure \ref{grid}) with Manhattan distances \citep{Ulmer2018} is assigned to a nurse (our method is able to handle any asymmetric area as well). We call it circular because all border regions have the same distance from the depot (nurse's home). Each day, the nurse's route begins and ends at the depot. The depot is assumed to be at the vertex of a region. The area includes $L$ square regions of equal size, where the center of each is used as a proxy for the location of all patients in the region (not an overly restrictive assumption provided that the regions are small enough). An additional small fixed value is included in each patient's service time to compensate for the travel time between patients' real locations in the same region as well as the distance between patients' real locations and the center of each region. Travel time from location $l$ to location ${l}'$ is known and equal to ${{d}_{l{l}'}},\forall l{l}'\in \{1,...,L\}\cup \{0\}$, where location $0$ represents the depot. 

A service type is prescribed to each new patient at admission. It is assumed that there are $K$ service types each of which is distinguished by a care pattern ${{h}_{k}}$, a probability distribution for the number of visits, a deterministic service time per visit ${{e}_{k}}$, and a wait-time target ${{T}_{k}}$. A care pattern states the number of days between consecutive visits to the patient, which we assume to be a constant. A wait-time target represents the maximum number of days an accepted new patient can wait for the first visit.

New referrals arrive dynamically each day following a known and steady-state probability distribution for each region and service type. Each nurse has a limited daily shift length of $\chi$ time units. The locations assigned to the nurse must, at a minimum, allow the nurse to travel to and from any of these locations and still have time to perform the longest visit in a single shift, i.e., ${{d}_{0,l}}+\underset{k}{\mathop{\max }}\,\{{{e}_{k}}\}+{{d}_{l,0}}\le \chi ,\forall l$.  

It is assumed that the number of visits of type $k$ is at most ${{J}_{k}}$. While the exact number of visits required for each patient is not known in advance, the status of the immediate next visit is determined at the time of the previous one {\textemdash} either at least one more visit is needed or there is no need for further visits. 

Each new referral can be accepted or rejected. The rejection decision can only be made for new referrals. In other words, no patient can be rejected in the middle of their treatment. A scheduled visit for each patient in a day can be done at any time during that day.

Each day the task of the home health care provider includes (i) determining whether to reject or accept each new referral, (ii) assigning the first visit of each accepted referral to a day within the wait-time target, (iii) deciding whether to divert any of today's visits, and finally (iv) finding the optimal route for serving today's visits. The cost function is a weighted sum of the rejection, diversion, overtime, and travel time costs.

\section{MDP formulation\label{sec:MDP_formulation}}
In this section, we formulate the problem as a discounted infinite-horizon \ac{MDP}, and establish Bellman's equations for it. 

\subsection{Decision epochs\label{sec:epoch}}
It is assumed that no new referral arriving in the middle of a day needs a visit in the same day. Therefore, decisions can be made once a day, at a time before the start of each day that is called a decision epoch. We assume a rolling planning horizon of $T$ days, where $T=\max \{\underset{k}{\mathop{\max }}\,\{{{h}_{k}}\},\underset{k}{\mathop{\max }}\,\{{{T}_{k}}\}\}$. In other words, day $t$ at the current decision epoch becomes day $t-1$ at the next decision epoch. This value of $T$ enables the model to keep track of the next visit for each patient within the planning horizon.

\subsection{State space\label{sec:state}}
At each decision epoch, the decision maker observes the patients assigned to each day in the planning horizon and the list of new referrals. A state is represented by $s=(\vec{x},\vec{y})$, where $\vec{x}=\{{{x}_{tklj}};t=1,...,T,k=1,...,K,l=1,...,L,j=0,...,{{J}_{k}}-1\}$ and $\vec{y}=\{{{y}_{kl}};\forall kl\}$. In each state, ${{x}_{tklj}}$ is the number of patients of type $k$ in location $l$ with $j$ visits completed before day $t$ and $(j+1)$-th visit already assigned to day $t$. We do not consider $j={{J}_{k}}$ because a patient with ${{J}_{k}}$ visits completed has already left the system. Also, ${{y}_{kl}}$ represents the number of new referrals of type $k$ from location $l$ waiting for a decision. 
An upper bound is defined for each state variable to make the state space finite \citep{Puterman1994,Adelman2004}. Therefore, the state space $S$ is defined as
\begin{equation}
	S=\{(\vec{x},\vec{y})|{{x}_{tklj}}\le x^{max}_{tklj},\forall tklj,{{y}_{kl}}\le y^{max}_{kl},\forall kl,{{x}_{tklj}}\in \mathbb{Z}_{\ge 0},\forall tklj,{{y}_{kl}}\in \mathbb{Z}_{\ge 0},\forall kl\},\label{eq:state}
\end{equation}
where $x^{max}_{tklj}$ and $y^{max}_{kl}$ represent the corresponding upper bound values for ${{x}_{tklj}}$ and ${{y}_{kl}}$, respectively. 

Due to the uncertain number of visits for each patient, only the immediate next visit of a patient is scheduled in the planning horizon. After visiting a patient of type $k$, the next visit, if needed, is booked on day ${{h}_{k}}$ of the planning horizon. Therefore, we do not have any visits of type $k$ booked on days later than $h_k$ in the next planning horizon for any patient who has been visited at least once, i.e., ${{x}_{tklj}}=0,\forall tklj,t\ge {{h}_{k}}+1,j\ge 1$. Furthermore, the first visit of a patient can be assigned to any day in the planning horizon within the wait-time target, i.e., ${{x}_{tkl,0}},\forall tkl,t\le {{T}_{k}}-1$ can be non-zero (because the planning horizon is rolling, we have ${{x}_{{{T}_{k}},kl,0}}=0,\forall kl$). For simplicity, we let the set of indexes for which ${{x}_{tklj}}$ is always zero be $\Theta =\{\forall tklj|(t\ge {{h}_{k}}+1,j\ge 1)\cup ({{T}_{k}}\le t\le T,j=0)\}$.

\subsection{Action set\label{sec:action}}
At each decision epoch, the scheduler takes an action $a=(\vec{r},\vec{n},\vec{z},\vec{o},u)$, where $\vec{r}=\{{{r}_{kl}};\forall kl\}$ denotes the number of new referrals of type $k$ from location $l$ that are rejected, $\vec{n}=\{{{n}_{tkl}};\forall tkl\}$ represents the number of accepted referrals of type $k$ from location $l$ assigned to day $t$, $\vec{z}=\{{{z}_{kl}};\forall kl\}$ gives the total number of visits of type $k$ from location $l$ that are diverted, $\vec{o}=\{{{o}_{l{l}'}};\forall l{l}'\}$ is a binary variable that takes the value of one if the nurse serves location $l$ before ${l}'$, and $u$ indicates the amount of overtime the nurse should work. The last three action variables,  $\vec{z},\vec{o}$, and $u$, affect only the first day of the planning horizon, and the decision to divert is only applied for one visit of each patient at a time. Any such action must satisfy the following constraints. According to constraint set \eqref{eq:r}, the number of rejections of each type for each location is equal to the number of new referrals minus the number of accepted referrals:
\begin{equation}
	{{r}_{kl}}={{y}_{kl}}-\sum\nolimits_{t=1}^{{{T}_{k}}}{{{n}_{tkl}},\forall kl}.\label{eq:r}
\end{equation}

Based on constraint set \eqref{eq:divert}, the number of diverted visits of each type for each location must be less than or equal to the corresponding number of visits assigned to day 1:
\begin{equation}
	{{z}_{kl}}\le\sum\nolimits_{j}{{{x}_{1,klj}}}+{{n}_{1,kl}},\forall kl.\label{eq:divert}
\end{equation}

To determine the optimal route to visit patients on the current day, we need to know the locations to be visited by the nurse on day 1. Through constraint sets \eqref{eq:w} and \eqref{eq:w_sign}, we set ${{w}_{l}}=1$ if there is at least one patient in location $l$ that must be visited by the nurse on day 1, where $M$ is a big number:
\begin{gather}
	\sum\nolimits_{kj}{{{x}_{1,klj}}}+\sum\nolimits_{k}{{{n}_{1,kl}}}-\sum\nolimits_{k}{{{z}_{kl}}}\le M{{w}_{l}},\forall l,\label{eq:w}
	\\
	{{w}_{l}}\in \{0,1\},\forall l.\label{eq:w_sign}
\end{gather}

Constraint set \eqref{f1} ensures that the total travel time from the depot to each location $l$ (captured by ${{f}_{l}}$) is at least equal to the direct travel time from the depot to that location: 
\begin{equation}
	{{f}_{l}}\ge {{d}_{0,l}},\forall l.\label{f1} \\
\end{equation}

According to constraint sets \eqref{f2} and \eqref{o_sign}, the total travel time from the depot to each location is greater than or equal to that of its preceding locations in the nurse’s route plus the travel time between the corresponding locations \citep{Sawik2000}:
\begin{gather}
	{{f}_{l}}+M(2+{{o}_{l{l}'}}-{{w}_{l}}-{{w}_{{{l}'}}})\ge {{f}_{{{l}'}}}+{{d}_{{l}'l}},\forall l{l}',l\ne {l}',\label{f2} \\
	{{o}_{l{l}'}}\in \{0,1\},\forall l{l}'.\label{o_sign}
\end{gather}

The nurse’s total travel time on day 1, shown by $q$, is obtained from constraint set \eqref{eq:total_travel} guaranteeing the return of the nurse to the depot at the end of the route with the addition of ${{d}_{l,0}}$:
\begin{equation}
	{{f}_{l}}+{{d}_{l,0}}\le q+M(1-{{w}_{l}}),\forall l.\label{eq:total_travel}
\end{equation}

Constraint \eqref{eq:g} determines the tour length of the nurse on day 1 denoted by $g$ as the total travel time plus the total service time of all patients visited by the nurse on day 1:
\begin{equation}
	g=q+\sum\nolimits_{k}{(\sum\nolimits_{lj}{{{x}_{1,klj}}}+\sum\nolimits_{l}{{{n}_{1,kl}}}-\sum\nolimits_{l}{{{z}_{kl}}}){{e}_{k}}}.\label{eq:g}
\end{equation}

Constraint \eqref{u1} ensures that $g$ does not exceed the regular shift length of the nurse (shown by $\chi $) plus overtime (represented by $u$), and finally, constraint \eqref{u2} imposes a maximum allowable overtime ${\chi }'$:
\begin{gather}
	g\le \chi +u,\label{u1} \\
	u\le {\chi }'.\label{u2}
\end{gather}

\subsection{State transition\label{sec:state_transition}}
The transition from state $s$ to state $s'=(\vec{{x}'},\vec{{y}'})$ (the state right before the start of the next day) depends on the chosen action for the current day and two stochastic elements --- the number of new referrals and the number of patients visited in the current day who need at least one more visit.

In the current day, the nurse visits ${{x}_{1,kl,0}}+{{n}_{1,kl}}$ patients of type $k$ in location $l$ for their first visit. In addition, she/he visits ${{x}_{1,klj}},\forall klj,j\ge 1$ patients for a follow-up visit leading to the completion of their ($j+1$)-th visit by the end of the day. We let the probability of needing at least one more visit for a patient of type $k$ after $j$ visits be equal to ${{p}_{k,j+1}}$. Mathematically, ${{p}_{kj}}=\PX({{\hat{J}}_{k}}\ge j|{{\hat{J}}_{k}}\ge j-1),\forall j\ge 1$, where ${{\hat{J}}_{k}}$ represents the random variable of the number of visits for service type $k$. Let $B(x,p)$ denote a binomial random variable with $x$ being the number of trials and $p$ the success probability. Therefore, $B({{x}_{1,kl,0}}+{{n}_{1,kl}},{{p}_{k,2}})$ and $B({{x}_{1,kl,j-1}},{{p}_{k,j+1}})$ respectively represent the corresponding random variable for the number of patients with one visit and $j\ge 2$ visits completed who are served on the current day and require at least one more visit. We denote the random variable corresponding to $\vec{{x}'}$ by $\vec{X}'$, which can be defined as
\begin{equation}
        \small
	X'_{tklj}(s,a)=\left\{ \begin{array}{*{35}{l}}
		B({{x}_{1,kl,0} }+{{n}_{1,kl} },{{p}_{k,2}}), & j=1,t={{h}_{k}}  \\
		B({{x}_{1,kl,j-1} },{{p}_{k,j+1}}), & j\ge 2,t={{h}_{k}}  \\
		{{x}_{t+1,klj} }, & j\ne 0,t\le {{h}_{k}}-1  \\
		{{x}_{t+1,kl,0} }+{{n}_{t+1,kl} }, & j=0,t<{{T}_{k}}  \\
		0, & \text{otherwise }(\text{i.e., } tklj\in \Theta )  \\
	\end{array} \right.,\forall s\in S,\forall a\in {{A}_{s}},\forall tklj\label{X_p},
\end{equation}
where ${{A}_{s}}$ represents the set of all feasible actions for a given state $s$, satisfying constraints \eqref{eq:state}-\eqref{u2}. Because ${{x}_{t+1,klj}}=0,\forall t=h_k,kl,j\ge 1$, the formulas for $j=1,t={{h}_{k}}$ and $j\ge 2,t={{h}_{k}}$ in equation \eqref{X_p} do not include ${{x}_{t+1,klj}}$. 

Consider $\vec{\xi }=\{{{\xi }_{klj}};\forall klj\}$ as the number of patients visited on the previous day and needing at least one more visit. Because the health statuses of different patients do not affect each other, it is logical to assume different elements of $\vec{\xi }$ are independent of each other. Similarly, the number of arrivals of different service types in the same region or in different regions of an area can be assumed independent of each other. Moreover, it can be assumed that the elements of $\vec{\xi }$ and $\vec{{y}'}$ are independent. Therefore, we have $\PX(s'|s,a)=\PX(\vec{\xi})\PX(\vec{{y}'})=\prod\limits_{klj}{\PX({{\xi }_{klj}})}\times \prod\limits_{kl}{\PX(y'_{kl})}$. 

\subsection{Cost function}
The immediate cost associated with a given state-action pair, $c(s,a)$, includes the costs of the rejection of new referrals, the diversion of visits of referrals already accepted, the nurse’s overtime work, and the travel time:
\begin{equation}
	c(s,a)=\sum\limits_{kl}{{{R}_{k}}{{r}_{kl} }}+\sum\limits_{kl}{{{Z}_{k}}{{z}_{kl} }}+Uu +Qq(s,a),\forall s\in S,a\in {{A}_{s}}\label{cost},
\end{equation}
where ${{R}_{k}}$ represents the cost of rejecting a referral of type $k$, ${{Z}_{k}}$ is the cost of diverting a visit of type $k$, $U$ is the unit cost of working overtime, and $Q$ denotes the unit cost associated with travel time. It is assumed that ${{R}_{k}}={{\zeta }^{r}}\EX[{{\hat{J}}_{k}}]{{e}_{k}}$, where  ${{\zeta }^{r}}$ is a weight factor, and $\EX[{{\hat{J}}_{k}}]$ represents the expected number of visits for a patient of type $k$. Thus, the rejection cost for each new referral of a given type is assumed to be proportional to the expected number of visits and the service time of that type. It represents the expected lost profit associated with not serving all the visits of a referral. 

For the diversion cost, we have ${{Z}_{k}}={{\zeta }^{z}}{{e}_{k}}$, where ${{\zeta }^{z}}$ is the corresponding weight (proportionally) factor. In a single-nurse setting, a diverted visit is assigned to an on-call nurse. Assuming that an on-call nurse is always available near each location, the diversion cost of a visit does not depend on its location. The diversion cost is two-fold --- the additional cost of serving a visit by an on-call nurse and the cost of violating continuity of care (because a different nurse serves the patient). Continuity of care increases the overall patient's quality of care and satisfaction \citep{Haggerty20031219}. 

In the cost function, the values of $U$ and $Q$ are assumed to be constant and respectively equal to ${{\zeta }^{u}}$ and ${{\zeta }^{q}}$. The weight values in the cost function (i.e., ${{\zeta }^{r}}$, ${{\zeta }^{z}}$, ${{\zeta }^{u}}$, and ${{\zeta }^{q}}$) are determined based on stakeholders’ preferences. We discuss this further in Section \ref{sec:general_settings}.
\subsection{Bellman's equations}\label{sec:Bel}
A policy $\mu$ maps each state $s \in S$ to an action $a \in A_s$. The value function of a given state $s$ under a policy $\mu$, denoted by $v_\mu(s)$, is defined as $v_\mu(s)=\EX[\sum_{i=0}^{\infty} \gamma^ic(s_i,\mu(s_i))|s_0=s]$, where $\gamma $ denotes the daily discount factor. A policy $\mu$ is defined to be better than or equal to a policy $\mu'$ if and only if $v_\mu(s) \le v_{\mu'}(s), \forall s \in S$. The corresponding value function of the optimal policy, shown by $v^*(s)$, is defined as $v^*(s)=\displaystyle \min_\mu v_\mu(s),\forall s\in S$. According to Bellman's optimality equations, $v^*(s)$ can be obtained by
\begin{equation}
	v^*(s)=\underset{a\in {{A}_{s}}}{\mathop{\min }}\,\{c(s,a)+\gamma \sum\limits_{{s}'\in S}{\PX({s}'|s,a)v^*({s}')}\},\forall s\in S\label{Bel}.
\end{equation}
Then, the optimal policy $\mu^*$ is obtained by {\small $\mu^*(s)=\underset{a \in A_s}{\mathop{\argmin }}\,\{c(s,a)+\gamma \sum\limits_{s' \in S}{\PX(s'|s,a)v^*({s}')}\},\forall s\in S$}.
\section{ALP approach}\label{sec:Approximate-linear-program}
In this section, we first transform Bellman’s equations into an equivalent \ac{LP}, create an ALP, and demonstrate how to optimize the ALP parameters in Section \ref{subsec:Obtaining-the-ALP}. Then, we formulate the action generation for a given state based on the resulting ALP parameters in Section \ref{subsec:ALP-based-action-generation}.
\subsection{ALP formulation}\label{subsec:Obtaining-the-ALP}
According to \ac{MDP} theory (see, e.g., \citealp{Puterman1994}, Section 6.9), solving Bellman's equations, given in \eqref{Bel}, is equivalent to solving the following LP model for any strictly positive set of state-relevance weights $\alpha (s),\forall s\in S$:
\begin{equation}
	\begin{split}
		& \mathbf{LP}:\underset{{\vec{v}}}{\mathop{\max }}\,\sum\limits_{s\in S}{\alpha (s)v(s)}, \\
		& \text{subject to} \\
		& c(s,a)+\gamma \sum\limits_{{s}'\in S}{\PX({s}'|s,a)}v({s}')\ge v(s),\forall s\in S, \forall a\in {{A}_{s}}.\\
	\end{split}\label{LP}
\end{equation}
Because there is a variable for every state and a constraint for every state-action pair, solving the LP is no easier than solving Bellman's equations. To tackle this issue, we approximate the value function by the following affine function of state variables:
\begin{equation}
	\tilde{v}(s)=\eta +\sum\limits_{tklj}{{{\tau }_{tklj}}{{x}_{tklj}} }+\sum\limits_{kl}{{{\rho }_{kl}}{{y}_{kl}} },\forall s\in S\label{eq:affine_fun},
\end{equation}
where ${{\tau }_{tklj}}\in \mathbb{R}_{\ge 0},\forall tklj$, ${{\rho }_{kl}}\in \mathbb{R}_{\ge 0},\forall kl$, and $\eta \in \mathbb{R}$ represent the \enquote{ALP parameters}. The optimal values of ${{\tau }_{tklj}}$ and ${{\rho }_{kl}}$ can be interpreted as the marginal expected infinite-horizon discounted cost of an additional visit and a new referral with the same indexes, respectively.  Let’s define $\vec{\tau }=\{{{\tau }_{tklj}};\forall tklj\}$ and $\vec{\rho }=\{{{\rho }_{kl}};\forall kl\}$. After utilizing $\tilde{v}(s)$ instead of $v(s)$ in the LP model, expansion, and rearrangement, we have the following primal model, called ALP:
\begin{subequations}
        \small
	\begin{align}
		& \textbf{ALP}:\underset{\eta ,\vec{\tau },\vec{\rho }}{\mathop{\max }}\,\left\{ \eta +\sum\limits_{tklj}{{{\tau }_{tklj}}\EX_{\alpha}[{{X}_{tklj}}]}+\sum\limits_{kl}{{{\rho }_{kl}}\EX_{\alpha}[{{Y}_{kl}}]} \right\}, \label{P-ALP-1} \\
		& \text{subject to} \nonumber \\
		& (1-\gamma )\eta +\sum\limits_{tklj}{{{\tau }_{tklj}}\left( {{x}_{tklj} }-\gamma \EX[X'_{tklj}(s,a)] \right)}+\sum\limits_{kl}{{{\rho }_{kl}}\left( {{y}_{kl} }-\gamma \EX[{{Y}_{kl}}] \right)}\le c(s,a),\forall s\in S,\forall a\in {{A}_{s}} \label{P-ALP-2}, \\
		& {{\tau }_{tklj}}\ge 0,\forall tklj;{{\rho }_{kl}}\ge 0,\forall kl, \label{P-ALP-3}
	\end{align}
\end{subequations}
where $\EX_{\alpha}[{{X}_{tklj}}]=\sum\limits_{s\in S}{\alpha (s){{x}_{tklj} }}$, $\EX_{\alpha}[{{Y}_{kl}}]=\sum\limits_{s\in S}{\alpha (s){{y}_{kl} }}$, $\EX[{{Y}_{kl}}]$ is the expected number of referrals of type $k$ from location $l$, and
\vspace{-\baselineskip} 
\begin{equation}
	\EX[X'_{tklj}(s,a)]=\left\{ \begin{array}{*{35}{l}}
		({{x}_{1,kl,0} }+{{n}_{1,kl} }){{p}_{k,2}}, & j=1,t={{h}_{k}}  \\
		{{x}_{1,kl,j-1} }{{p}_{k,j+1}}, & j\ge 2,t={{h}_{k}}  \\
		{{x}_{t+1,klj} }, & j\ne 0,t\le {{h}_{k}}-1  \\
		{{x}_{t+1,kl,0} }+{{n}_{t+1,kl} }, & j=0,t<{{T}_{k}}  \\
		0, & \text{otherwise}  \\
	\end{array} \right.,\forall s\in S,\forall a\in {{A}_{s}},\forall tklj.\label{E_X_p}
\end{equation}

Due to the large number of constraints in the ALP for any realistically-sized problem, we derive its dual, called D-ALP, and solve it using column generation. The D-ALP can be written as:
\begin{subequations}
	\begin{align}
		& \textbf{D-ALP}:\underset{{\vec{\beta }}}{\mathop{\min }}\,\sum\limits_{s\in S,a\in {{A}_{s}}}{\beta (s,a)c(s,a)}, \\
		& \text{subject to}\nonumber \\
		& (1-\gamma )\sum\limits_{s\in S,a\in {{A}_{s}}}{\beta (s,a)}=1,\label{D_ALP_1} \\
		& \sum\limits_{s\in S,a\in {{A}_{s}}}{\beta (s,a)({{x}_{tklj} }-\gamma \EX[X'_{tklj}(s,a)])}\ge \EX_{\alpha}[{{X}_{tklj}}],\forall tklj\label{D_ALP_2}, \\
		& \sum\limits_{s\in S,a\in {{A}_{s}}}{\beta (s,a)({{y}_{kl} }-\gamma \EX[{{Y}_{kl}}])}\ge \EX_{\alpha}[{{Y}_{kl}}],\forall kl,\label{D_ALP_3} \\
		& \vec{\beta }\ge 0,
	\end{align}
\end{subequations}
where $\beta (s,a)$ represents the dual variable for each state-action pair and $\vec{\beta }=\{\beta (s,a);\forall s\in S,\forall a\in {{A}_{s}}\}$. The column generation method starts by generating a feasible state-action pair and adding it to a set $S'$. Then, it solves the D-ALP considering only the state-action pair in $S'$. Next, using the dual prices of \eqref{D_ALP_1}, \eqref{D_ALP_2}, and \eqref{D_ALP_3} as estimates for $\eta, {{\tau }_{tklj}},$ and ${{\rho }_{kl}}$, it finds the most violated constraint in the primal ALP, and adds the state-action pair associated with this constraint into the set $S'$ before re-solving the D-ALP with the updated $S'$. The process iterates until no primal constraint is violated.
To obtain the state-action pair corresponding to the most violated constraint at each iteration of the column generation method, the following sub-problem, called S-ALP, must be solved:
\begin{equation}
    \footnotesize
    \textbf{S-ALP}:\underset{s\in S,a\in {{A}_{s}}}{\mathop{\min }}\,\left[ c(s,a)-\eta (1-\gamma )-\sum\limits_{tklj}{{{\tau }_{tklj}}\left( {{x}_{tklj} }-\gamma \EX[{{{{X}'}}_{tklj}(s,a)}] \right)}-\sum\limits_{kl}{{{\rho }_{kl}}\left( {{y}_{kl} }-\gamma \EX[{{Y}_{kl}}] \right)} \right].\label{S_ALP}
\end{equation}
If \eqref{S_ALP} is negative, then the resulting state-action pair is added to $S'$, and the next iteration is performed; otherwise, the current estimates of the ALP parameters are optimal. Unlike the original (non-approximated) LP, the solution to the ALP may be sensitive to the choice of $\alpha$ \citep{DeFarias2003850, Adelman2004}. For simplicity, we assume $\EX_{\alpha}[{{X}_{tklj}}]=\varepsilon ,\forall tklj\notin \Theta $, and zero otherwise, where $\varepsilon $ is an unknown constant. We provide a procedure to tune the choice of $\varepsilon $ in Section \ref{sec:ALP_tuning}. It should be mentioned that ${{E}_{\alpha}}[{{Y}_{kl}}],\forall kl$ is not a sensitive parameter, and we set it equal to $\EX[{{Y}_{kl}}]$.

To obtain the initial state-action pair required to launch column generation, we let $\beta (s,a)=1/(1-\gamma )$ based on \eqref{D_ALP_1}. Then, we assume that all referrals are accepted (i.e., ${{r}_{kl}}=0,\forall kl$ and ${{y}_{kl}}=\sum\limits_{t}{{{n}_{tkl}}},\forall kl$) and all visits assigned to day 1 are diverted (i.e., ${{z}_{kl}}=\sum\limits_{j}{{{x}_{1,klj}}}+{{n}_{1,kl}},\forall kl$ and $u=q=0$). Any solution satisfying constraints \eqref{D_ALP_2} and \eqref{D_ALP_3} considering the assumptions mentioned above provides an initial state-action pair.

With this method, we are able to obtain the optimal ALP parameters for small-scale instances of the problem but solving large-scale ones may be very time-consuming or even intractable. To illustrate this issue, the execution time to derive the optimal ALP parameters of an instance with 6 regions in a rectangular area ($L=6$), 1 service type ($K=1$), and 3 visits per referral on average ($\EX[\hat{J}_1]=3$) may be around one hour while that of each of the instances ($L=24,K=1,\EX[\hat{J}_1]=3$), ($L=6,K=1,\EX[\hat{J}_1]=12$), and ($L=6,K=2,\EX[\hat{J}_k]=3,\forall k$) may be more than a day. We present two heuristic reduction techniques to tackle this issue in Section \ref{sec:Acceleration-techniques-for}.
\subsection{ALP policy}\label{subsec:ALP-based-action-generation}
Having a set of ALP parameters ${{\pi }^{*}}=({{\eta }^{*}},\tau _{tklj}^{*},\rho _{kl}^{*})$, the ALP policy obtains an action for any given state by solving the approximate version of equation \eqref{Bel} that can be re-written as
\begin{equation}
	\underset{a\in {{A}_{s}}}{\mathop{\min }}\,\left\{ c(s,a)+\gamma \left( \eta +\sum\limits_{tklj}{\tau _{tklj}^{*}\EX[{{{{X}'}}_{tklj}}]}+\sum\limits_{kl}{\rho _{kl}^{*}\EX[{{Y}_{kl}}]} \right) \right\}.\label{AG1}
\end{equation}
After expansion and integrating the constant values, we have:
\vspace{-\baselineskip} 
\begin{multline}
	\underset{a\in {{A}_{s}}}{\mathop{\min }}\,\left\{ \sum\limits_{kl}{{{R}_{k}}{{r}_{kl}}}+\sum\limits_{kl}{{{Z}_{k}}{{z}_{kl}}}+Uu+Qq+ \right. \\
	\left. \gamma \left( \sum\limits_{tkl;t={{h}_{k}}}{\tau _{tkl,1}^{*}{{p}_{k,2}}{{n}_{1,kl}}}+\sum\limits_{tkl;t\le {{T}_{k}}-1}{\tau _{tkl,0}^{*}{{n}_{t+1,kl}}} \right) \right\}+\text{constant}.\label{AG2}
\end{multline}
As a result of \eqref{AG2}, Proposition \ref{Prop_radius} will help us determine regions for which the ALP policy for a given set of the ALP parameters always accepts or rejects new referrals of a given service type.
\begin{proposition}
	Any referral of type $k$ coming from location $l$ is always rejected by the ALP policy if $\min \{\gamma \tau _{{{h}_{k}},kl,1}^{*}{{p}_{k,2}},\underset{t,t\le {{T}_{k}}-1}{\mathop{\min }}\,\{\gamma \tau _{tkl,0}^{*}\}\}>{{R}_{k}}$, and is always accepted if $\underset{t,t\le {{T}_{k}}-1}{\mathop{\min }}\,\{\gamma \tau _{tkl,0}^{*}\}<{{R}_{k}}$.\label{Prop_radius}
\end{proposition}
\noindent\textsc{Proof}. The three possible decisions for any referral of type $k$ from location $l$ are: (i) rejection, (ii) acceptance and assignment to day 1, (iii) acceptance and assignment to day $t,2\le t\le {{T}_{k}}$. The additional costs incurred by (i) and (iii) in \eqref{AG2} are known and equal to ${{R}_{k}}$ and $\underset{t,t\le {{T}_{k}}-1}{\mathop{\min }}\,\{\gamma \tau _{tkl,0}^{*}\}$ in this order. Moreover, the resulting additional cost of (ii) is at least equal to $\gamma \tau _{{{h}_{k}},kl,1}^{*}{{p}_{k,2}}$. Therefore, if $\underset{t,t\le {{T}_{k}}-1}{\mathop{\min }}\,\{\gamma \tau _{tkl,0}^{*}\}< {{R}_{k}}$ holds, we are sure that assignment to day $t,2\le t\le {{T}_{k}}$ would be better than rejection, and we would thus never reject. On the other hand, if the minimum possible additional cost of assignment, i.e., $\min \{\gamma \tau _{{{h}_{k}},kl,1}^{*}{{p}_{k,2}},\underset{t,t\le {{T}_{k}}-1}{\mathop{\min }}\,\{\gamma \tau _{tkl,0}^{*}\}\}$ is greater than ${{R}_{k}}$, then we never accept them.

\section{A closed-form solution for a special case\label{sec:A-closed-form-solution}}
Consider a special case of the problem, called ${\cal P}_s$, in which all referrals are accepted, no overtime is allowed, the diversion cost is much greater than the travel cost, i.e., ${{\zeta }^{z}} \gg {{\zeta }^{q}}$, travel times are symmetric for each pair of locations, i.e., ${{d}_{l{l}'}}={{d}_{{l}'l}},\forall l{l}'$, and service times are the same for all service types, i.e., ${{e}_{k}}=e,\forall k$. Theorem \ref{theorem1} presents a closed-form solution for the optimal ALP parameters of ${\cal P}_s$ (denoted by ${{\eta }^{*}}({\cal P}_s)$, $\tau _{tklj}^{*}({\cal P}_s)$, and $\rho _{kl}^{*}({\cal P}_s)$) when
\begin{equation}
	\sum\nolimits_{klj}{\beta _{1,klj}^{x} +\sum\nolimits_{kl}{\beta _{1,kl}^{n} }<\frac{1}{1-\gamma }}\label{T1_cond}
\end{equation}
is satisfied, where $\beta _{1,klj}^{x}=\sum\limits_{(s,a)\in \Lambda ,{{x}_{1,klj}}>0}{\beta^* (s,a)}$,  $\beta _{1,kl}^{n}=\sum\limits_{(s,a)\in \Lambda ,{{n}_{1,kl}}>0}{\beta^* (s,a)}$, $\Lambda$ is a set of state-action pairs for which the primal solution (${{\eta }^{*}}({\cal P}_s)$, $\tau _{tklj}^{*}({\cal P}_s), {\forall tklj}$, $\rho _{kl}^{*}({\cal P}_s), \forall {kl}$) generates binding constraints in the ALP model, and $\beta^* (s,a), \forall a\in {{A}_{s}},\forall s\in S$ is the corresponding dual solution for this primal one. More details about how to calculate $\beta _{1,klj}^{x}$ and $\beta _{1,kl}^{n}$ together with a proof for Theorem \ref{theorem1} are provided in Section \ref{sec:proof-theorem1}.

\begin{theorem}\label{theorem1}
When condition \eqref{T1_cond} holds, the optimal ALP parameters of ${\cal P}_s$ are given by
\end{theorem}
\vspace{-\baselineskip} 
\begin{subequations}
	\begin{align}
		& \tau _{tklj}^{*}({\cal P}_s)=\left\{ \begin{array}{*{35}{l}}
			{{\gamma }^{t-1}}\tau _{1,kl,{{J}_{k}}-1}^{*}({\cal P}_s)\sum\limits_{i=0}^{{{J}_{k}}-1-j}{({{\gamma }^{i{{h}_{k}}}}\underset{{i}'=0}{\overset{i-1}{\mathop{\prod }}}\,{{p}_{k,j+2+{i}'}})}, & \forall tklj\notin \Theta   \\
			0, & \text{otherwise}  \\
		\end{array}, \right.\label{T1_1} \\
		& \tau _{1,kl,{{J}_{k}}-1}^{*}({\cal P}_s)={{\zeta }^{q}}\frac{2{{d}_{0,l}}}{\left\lfloor (\chi -2{{d}_{0,l}})/e \right\rfloor },\forall kl,\label{T1_2} \\
		& \rho _{kl}^{*}({\cal P}_s)={{\gamma }^{{{T}_{k}}-1}}\tau _{1,kl,{{J}_{k}}-1}^{*}({\cal P}_s)\sum\limits_{i=0}^{{{J}_{k}}-1}{({{\gamma }^{i{{h}_{k}}}}\underset{{i}'=0}{\overset{i-1}{\mathop{\prod }}}\,{{p}_{k,2+{i}'}})},\forall kl,\label{T1_3} \\
		& {{\eta }^{*}}({\cal P}_s)=\frac{\gamma }{1-\gamma }\sum\nolimits_{kl}{\rho _{kl}^{*}({\cal P}_s)\EX[{{Y}_{kl}}]}.\label{T1_4}
	\end{align}
\end{subequations}

This closed-form solution can be explained intuitively. Equation \eqref{T1_1} implies that $\tau _{tklj}^{*}({\cal P}_s)$ is equal to the expected present value of the marginal costs of the remaining visits of a patient of type $k$ from location $l$ with $j$ visits completed before day $t$ and $(j+1)$-th visit assigned to day $t$. The marginal costs of all visits of a patient of type $k$ from location $l$ are considered to be the same and equal to $\tau _{1,kl,{{J}_{k}}-1}^{*}({\cal P}_s)$, which corresponds to the marginal cost of a single visit (the last possible one) of the same type and location assigned to day 1. There are at most $J_k-j$ visits left for this patient that should be done on days $t$, $t+h_k$, ..., and $t+(J_k-1-j)h_k$, respectively. A patient with $j$ visits completed and $(j+1)$-th visit assigned to a day within the wait-time target would need an additional visit (i.e., the $(j+1)$-th one) with probability one, two additional visits with probability ${{p}_{k,j+2}}$, ..., and $J_k-j$ ones with probability $\underset{{i}'=0}{\overset{J_k-j-2}{\mathop{\prod }}}\,{{p}_{k,j+2+{i}'}}$.

The value of $\tau _{1,kl,{{J}_{k}}-1}^{*}({\cal P}_s)$ can be obtained from \eqref{T1_2}. The right-hand side can be interpreted as the minimum travel time cost per visit in location $l$, obtained when the maximum possible number of visits are served in this location. The term $2{{d}_{0,l}}$ indicates the minimum travel time when visiting at least one patient in location $l$, which is equal to the travel time to and from the depot to location $l$. The term $\left\lfloor (\chi -2{{d}_{0,l}})/e \right\rfloor $ represents the maximum number of patients that could be visited if the nurse only visited location $l$. This nicely links the travel cost to the marginal cost of a scheduled patient.

Equation \eqref{T1_3} indicates that $\rho _{kl}^{*}({\cal P}_s)$ is equal to the marginal discounted cost of a new referral of type $k$ from location $l$ whose first visit is assigned to the latest possible day, i.e., ${{T}_{k}}$. The right-hand side of this equation is obtained from \eqref{T1_1} when $t=T_k$ and $j=0$. Knowing the values of $\tau _{tklj}^{*}({\cal P}_s),\forall tklj$ and $\rho _{kl}^{*}({\cal P}_s),\forall kl$, and based on the proof presented in Section \ref{sec:proof-theorem1}, we can calculate ${{\eta }^{*}}({\cal P}_s)$ using \eqref{T1_4}.

According to \eqref{T1_1}, $\tau _{tklj}^{*}={{\gamma }^{t-1}}\tau _{1,klj}^{*},\forall tklj\notin \Theta $ holds for the problem ${\cal P}_s$ when condition \eqref{T1_cond} is satisfied. As a result of this property, Proposition \ref{Prop_assignment} characterizes the assignment day of referrals.
\begin{proposition}
	When $\tau _{tkl,0}^{*}={{\gamma }^{t-1}}\tau _{1,kl,0}^{*},\forall tkl\notin \Theta $ holds, the ALP policy assigns an accepted referral to the last day within the wait-time target if it is not assigned to day 1.\label{Prop_assignment}
\end{proposition}
\noindent\textsc{Proof}. Assignment of a referral to different days in the range $[2 - {{T}_{k}}]$ only affects the term $\sum\limits_{tkl;t\le {{T}_{k}}-1}{\tau _{tkl,0}^{*}{{n}_{t+1,kl}}}$ in \eqref{AG2}. Because $\tau _{tkl,0}^{*}={{\gamma }^{t-1}}\tau _{1,kl,0}^{*},\forall tkl\notin \Theta $ holds, we have $\underset{t,2\le t\le {{T}_{k}}}{\mathop{\arg \min }}\,\{\tau _{tkl,0}^{*}\}={{T}_{k}}$. 
\section{Heuristic reduction techniques for the ALP}\label{sec:Acceleration-techniques-for}
While determining the optimal form for a special case is interesting, it does not directly provide us with a solution to the more realistic, general problem (denoted by ${\cal P}_g$). In this section, we provide two heuristic reduction techniques that enable us to derive the ALP parameters for the general problem. These techniques were inspired by the closed-form solution of the special case given in Section \ref{sec:A-closed-form-solution}. We call these reduction techniques heuristics because we do not provide any proof of their performance. However, we demonstrate through the numerical results in Section \ref{sec:effects_of_accel} that their application does not result in any loss of optimality for the random instances we generated.
\subsection{A two-index $\tau$ instead of the original four-index one}\label{sec:two_index}
From \eqref{T1_1}, {\footnotesize $\tau_{tklj}^{*}({\cal P}_s)={{\delta}_{tjk}}\tau_{1,kl,{{J}_{k}}-1}^{*}({\cal P}_s),\forall tklj\notin \Theta $}, where {\footnotesize $\delta_{tjk} = {{\gamma}^{t-1}}\sum\limits_{i=0}^{{{J}_{k}}-1-j}{({{\gamma}^{i{{h}_{k}}}}\displaystyle \prod_{i'=0}^{i-1}{{p}_{k,j+2+i'}})},\forall tjk$} is a known constant. For the general problem, $\tau_{tklj}^{*}({\cal P}_g)={{\delta }_{tjk}}\tau_{1,kl,{{J}_{k}}-1}^{*}({\cal P}_g),\forall tklj\notin \Theta $ held for all small-scale instances we solved but $\tau _{1,kl,{{J}_{k}}-1}^{*}({\cal P}_g)$ was not necessarily equal to $\tau _{1,kl,{{J}_{k}}-1}^{*}({\cal P}_s)$. Therefore, as a heuristic reduction technique, we modify the ALP model, given in \eqref{P-ALP-1}-\eqref{P-ALP-3}, by replacing ${{\tau }_{tklj}}$ with ${{\delta }_{tjk}}{{{\tau }'}_{kl}}$. We call the resulting modified primal model and its dual model ALP-2I and D-ALP-2I, respectively. Then, the approximately optimal ALP parameters can be derived through the D-ALP-2I model with the sub-problem S-ALP. There are $TKL\displaystyle\sum_k J_k-KL$ fewer constraints in D-ALP-2I than D-ALP making this technique extremely helpful especially when either the $T$ or $J_k,\forall k$ values are large.
\subsection{Transforming the area into a one-dimensional array}\label{sec:1D}
The larger the number of regions, $L$, in a specific area assigned to a nurse, the more precise will be the travel times in our model but also the more computationally intensive. To solve instances with a large $L$, we present the following heuristic reduction technique. As seen in Section \ref{sec:A-closed-form-solution}, the optimal ALP parameters for different regions with the same distance from the depot are the same for ${\cal P}_s$. This property was satisfied for all small-scale instances we solved for ${\cal P}_g$ as well. After modifying the ALP model of ${\cal P}_g$ based on having the same ALP parameters for the regions with the same distance from the depot, the resulting model inspired a proxy problem for ${\cal P}_g$, with much fewer regions, that reproduced the optimal ALP parameters of this problem for all small-scale instances we checked (Section \ref{sec:effects_of_accel}).

\begin{figure}[h]
	\begin{center}
		\includegraphics[height=1.7in]{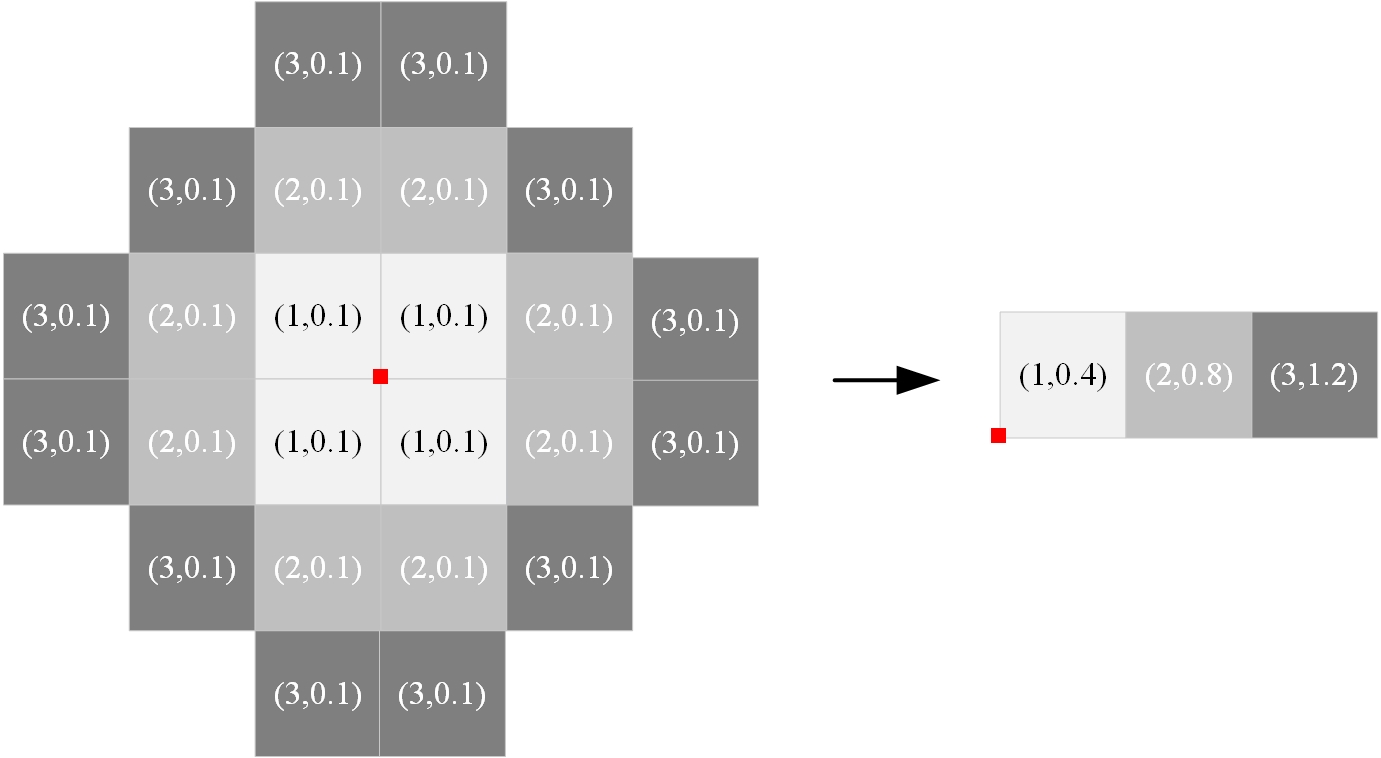}
		\caption{A \ac{2D} area and its corresponding \ac{1D} array} \label{fig-transform}
	\end{center}
\end{figure}

This proxy problem, called ${\cal P}^1_g$, includes a \ac{1D} array with ${{L}^{*}}$ regions produced by merging all regions with the same distance from the depot in the \ac{2D} area associated with ${\cal P}_g$. The \ac{1D} array consists of one row of regions starting from the depot. The distance of each region from the depot in the \ac{1D} area is equal to that of its corresponding regions in the \ac{2D} area.  
An example of transforming a \ac{2D} area into a \ac{1D} array can be seen in Figure \ref{fig-transform}.  We have assumed that there is only one service type in this example. The numbers in parentheses at the center of each region $l$ represent the distance from the depot and the referral rate of the only service type in that region. As an example, the distance from the depot to all of the regions with the lightest color in the \ac{2D} area is equal to one, and the arrival rate for the corresponding region in the \ac{1D} area is equal to the sum of the arrival rates of all of these regions. We can show that a 1D array with ${{L}^{*}}$ regions can represent a \ac{2D} area with at most $2{{L}^{*}}({{L}^{*}}+1)$ regions. 

Define ${\cal L}_{l^*},\forall {l^*}\in \{1,...,{{L}^{*}}\}$ as the set of all regions in the \ac{2D} area for which the distance of each to the depot is equal to ${{d}_{0,{l^*}}}({\cal P}^1_g)$, where ${{d}_{0,{l^*}}}({\cal P}^1_g)$ represents the distance of region ${l^*}$ from the depot in the \ac{1D} area. We have ${{\lambda }_{k{l^*}}}({\cal P}^1_g)=\sum\limits_{\forall l\in {{\cal L}_{{{{l^*}}}}}}{{{\lambda }_{kl}}}({\cal P}_g),\forall k{l^*}$, $\EX_{\alpha}[{{X}_{tk{l^*}j}}]({\cal P}^1_g)=\sum\limits_{l\in {{\cal L}_{{{{l^*}}}}}}{\EX_{\alpha}[{{X}_{tklj}}]}({\cal P}_g),\forall tkj{l^*}$, and 
$\EX_{\alpha}[{{Y}_{k{l^*}}}]({\cal P}^1_g)=\sum\limits_{l\in {{\cal L}_{{{{l^*}}}}}}{\EX_{\alpha}[{{Y}_{kl}}]}({\cal P}_g),\forall k{l^*}$.

The proxy ${\cal P}^1_g$ may be solved using the ALP or the ALP-2I to get the optimal (or approximately optimal) ALP parameters denoted by $\tilde{\eta }({\cal P}^1_g)$, ${{\tilde{\tau }}_{tklj}}({\cal P}^1_g)$, and ${{\tilde{\rho }}_{kl}}({\cal P}^1_g)$. Then, the approximately optimal ALP parameters for ${\cal P}_g$ are obtained by
$\tilde{\eta }({\cal P}_g)=\tilde{\eta }({\cal P}^1_g)$, 
${{\tilde{\tau }}_{tklj}}({\cal P}_g)={{\tilde{\tau }}_{tk{l^*}j}}({\cal P}^1_g),\forall tkjl,l\in {\cal L}_{l^*},\forall {l^*}$, and 
${{\tilde{\rho }}_{kl}}({\cal P}_g)={{\tilde{\rho }}_{k{l^*}}}({\cal P}^1_g),\forall kl,l\in {\cal L}_{l^*},\forall {l^*}$.
Deriving the ALP parameters of ${\cal P}_g$ using the instruction described above is labeled ALP-1D or ALP-1D-2I depending on the model (ALP or ALP-2I) used to obtain the ALP parameters of the corresponding ${\cal P}^1_g$. With this technique, we need only solve a much smaller primal model and sub-problem using column generation.

\section{Numerical results}\label{sec:Numerical-results}
In this section, we first describe the general settings for the numerical experiments in Section \ref{sec:general_settings}. Then, we evaluate the performance of the ALP policy in Section \ref{sec:ALPvsMyopic} and illustrate some features
of this policy in Section \ref{sec:manage}. 

We used CPLEX 12.8 to solve all the mathematical models in the paper. All algorithms were coded in C\#. The numerical results described below were obtained using a PC with an Intel(R) Core (TM) i9-10900X CPU running at 3.70 GHz and 64 GB of RAM.

The ALP parameters can be obtained by solving the original ALP given in \eqref{P-ALP-1}-\eqref{P-ALP-3} or one of its variants mentioned in Section \ref{sec:Acceleration-techniques-for}, i.e., ALP-2I, ALP-1D and ALP-1D-2I. Based on a numerical experiment reported in Section \ref{sec:effects_of_accel}, the ALP-1D-2I obtained the optimal ALP parameters in the shortest time for all tested small instances. Therefore, we use this variant of the ALP in all experiments given in this section. 

\subsection{General settings}\label{sec:general_settings}Table \ref{table-general-settings} outlines the control parameters together with their levels/ranges and default values. Choices were made to evaluate the performance of the ALP across a wide variety of home care organizational settings with varying types of patients and goals.

\newcolumntype{M}[1]{>{\centering\arraybackslash}m{#1}}
\begin{table}[h]
	\caption{Control parameters in the numerical experiments}
	\footnotesize
	\centering
	\resizebox{\textwidth}{!}{%
		\fontsize{7}{10}\selectfont
		\begin{tabular}{M{3.5cm}M{3.5cm}M{5.0cm}M{2.5cm}M{2.5cm}}
			\hline
			Scope                               & Parameters                   & Range/Levels                                                  & Default for $K=1$ & Default for $K>1$ \\ \hline
			\multirow{4}{*}{Area}               & Distance type                & Manhattan                                                     & Manhattan         & Manhattan         \\  
			& Area shape                   & Circular, Rectangular                                                       & Circular            & Circular            \\ 
			& Diameter of the area (hours) & \{0.5, 1, 2\}                                                   & 0.5               & 0.5               \\  
			& $L$                            & \{6, 24, 60\}                                                     & 24                & 24                \\ \hline
			Shift   length and overtime & $(\chi ,{\chi }')$ (hours)          & \{(8, 0), (8, 2)\}                                       & (8, 0)             & (8, 2)             \\ \hline
			\multirow{6}{*}{Service types} &
			${{h}_{k}}$ (days) &
			Fixed or random with ${{h}_{k}}\in   \{1,...,7\}$ &
			1 &
			Random \\ 
			&
			${{\hat{J}}_{k}}$ &
			Deterministic, truncated Poisson, or discrete uniform with mean 6, 8, or 10 &
			Truncated Poisson with mean 8 &
			Truncated Poisson with random mean \\ 
			& ${{e}_{k}}$ (hours)          & Fixed or random with ${{e}_{k}}\in   \{0.25, 0.50, 0.75, 1.00\}$ & 0.5               & Random            \\ 
			& ${{T}_{k}}$                  & Fixed or random with ${{T}_{k}}\in   \{1,...,5\}$             & 5                 & Random            \\ \hline
			\multirow{2}{*}{Referrals}          & $\bar{D}$ (hours)                       & \{7.5, 8.0, 8.5, 9.0, 9.5, 15.0\}                                             & 8.5               & 9.5               \\ 
			& ${{\lambda }_{kl}}$          & Fixed or random             & Fixed             & Random            \\ \hline
			\multirow{2}{*}{Cost   function} &
			$({{\zeta }^{r}},{{\zeta   }^{z}},{{\zeta }^{u}},{{\zeta }^{q}})$ &
			\{(5, 7.5, 2, 0.1), (5, 10, 2, 0.1),\space \space \space \space \space \space \space \space \space \space \space \space \space \space \space(5, 20, 2, 0.1), (5, 1000, 2, 0.1)\} &
			(5, 10, 2, 0.1) &
			(5, 10, 2, 0.1) \\  
			& $\gamma $                    & 0.99                                                          & 0.99              & 0.99              \\ \hline
		\end{tabular}\label{table-general-settings}%
	}
\end{table}
Service areas of 30 minutes in diameter in metropolitan cities generate sufficient volume for a single nurse \citep{Bennett2011}. For our base case scenario, we consider a circular area with Manhattan distances, a diameter of 30 minutes, and a depot at the center. We evaluate the sensitivity of the results to larger diameters of 1 and 2 hours. The default number of regions in the area (i.e., $L$) is 24 but we also report the results for an example with 60 regions in Section \ref{sec:manage}. Moreover, because of the long execution times for the model without the heuristic reduction techniques, we use a rectangular area with 6 regions to evaluate the heuristic reduction techniques for the ALP in Section \ref{sec:effects_of_accel}. We consider a regular shift length of 8 hours for a nurse (i.e., ${\chi}=8$), and a regular overtime limit of 2 hours (${\chi}'=2$) when it is allowable.

In our base case scenario, we assume that the number of visits for each service type follows a truncated Poisson distribution, and the minimum, maximum, and expected value of the number of visits are known. We let the minimum number of visits of any service type equal to one because we assume that any accepted referral would need at least one visit. Furthermore, we let the maximum number of visits of type $k$ be equal to $3{{\bar{J}}_{k}}$, where $\bar{J}_k$ represents the nearest integer to $\EX[\hat{J}_k]$. For the instance used to estimate the optimality gap of the ALP policy in Appendix B, it is exceptionally equal to $2{{\bar{J}}_{k}}$ due to tractability issues. In the home care literature, it is typical to plan for a duration ranging from 1 to 4 weeks (see, e.g., \citealp{Cappanera2018} and \citealp{Demirbilek2019a}). Given the various care patterns under consideration, we selected three values of 6, 8, and 10 for $\bar{J}_k$ in our experiments. Moreover, we consider two additional scenarios regarding the number of visits for each service type. In one scenario, it is assumed to be known, while in the other, it follows a discrete uniform distribution, with the maximum number of visits for type $k$ being $2\bar{J}_{k}-1$. Hence, we investigate how varying levels of variance in the distribution of the number of visits for each service type --- including \enquote{zero} (deterministic), \enquote{medium} (truncated Poisson), and \enquote{high} (discrete uniform) --- impact performance of the ALP policy.

In practice, most visits require a service time no longer than 1 hour \citep{Nickel2012}. We select the service time for each service type (i.e., $e_k$) among the values of 15, 30, 45, and 60 minutes. The default value of the service time for experiments with only one service type (i.e., when $K=1$) is 30 minutes \citep{Bennett2011,Demirbilek2019a}. We assume that the first visit of an accepted new referral of type $k$ must be done within ${{T}_{k}},\forall k$ days, where ${{T}_{k}}\in \{1,...,5\},\forall k$ depending on the urgency level. The maximum time (in days) between consecutive visits for a patient is assumed to be one week, i.e., ${{h}_{k}}\le 7,\forall k$. 

We define the average daily new demand ($\bar{D}$) as the expected value of the sum of the service times of all visits of all daily referrals, i.e., $\bar{D}=\sum\nolimits_{kl}{{{\lambda }_{kl}}}{{e}_{k}}{{\bar{J}}_{k}}$. In the base case scenario for instances with one service type, $\bar{D}$ is equal to 8.5 hours and overtime is not allowed (we analyze the effect of overtime in a separate experiment). As the nurse's tour length for each day is equal to the sum of the service times and travel times, this default value would create a sufficiently congested system. Moreover, we investigate $\bar{D}$ equal to 8.0 and 9.0 hours for two instances with one service type, $\bar{D}=9.5$ hours for instances with multiple service types, $\bar{D}=15.0$ hours in an instance with a larger area than the default one used in Section \ref{sec:manage}, and $\bar{D}=7.5$ hours in an instance used to estimate the optimality gap of the ALP policy in Appendix B.

In similar problem settings, the number of new arrivals is typically modeled as a Poisson process (see, e.g., \citealp{Patrick2008,Angelelli2009,Ulmer2018}, and \citealp{Demirbilek2019a}). We make the same assumption here and also assume that the daily referral rate for each service type and region, i.e., ${{\lambda }_{kl}},\forall kl$, does not change during the simulation period. We consider two scenarios for the arrival rates: (i) the same referral rate for all service types and regions, i.e., ${{\lambda }_{kl}}=\bar{\lambda },\forall kl$ ($\bar{\lambda }$ is a constant determined with respect to the $\bar{D}$); and (ii) a uniformly distributed random referral rate for each service type and region in the range $[0,2\bar{\lambda }]$. Because the referral rates are generated randomly in scenario (ii), we need to scale the resulting referral rates to meet the prespecified $\bar{D}$.

The daily discount factor is set equal to 0.99 \citep{Patrick2008}. The weights of the different components in the cost function \eqref{cost} are set at $({{\zeta }^{r}},{{\zeta }^{z}},{{\zeta }^{u}},{{\zeta }^{q}})=(5, 10, 2, 0.1)$ though we do a sensitivity analysis with regards to these values as well. It is important to make a number of justifications about the corresponding values. First, if we do not have ${{\zeta }^{z}}>{{\zeta }^{r}}$, then it would not be worth rejecting any new referral as the expected rejection cost would always be greater than the total diversion cost of the expected number visits. Because diversions violate continuity of care and impose additional costs, companies would rather avoid them as much as possible. Second, to reject the least possible number of referrals, doing even all visits of a patient in overtime should have more preference than rejecting her/him. For this purpose, we should have ${{\zeta }^{u}}<{{\zeta }^{r}}$. Third, the weight of travel time should be much smaller than the other weights to avoid any rejection of a referral or diversion of a visit because of travel time. We tend to serve any far location if it is possible considering the shift length limitation. The travel time matters when we have multiple best solutions in terms of the total of the other cost items.

\subsection{Evaluation of the ALP policy}\label{sec:ALPvsMyopic}
In this section, we compare the ALP policy with a myopic policy and a scenario-based policy (hereafter referred to as \enquote{Myopic} and \enquote{SB}). 

In our partner home care company, the decision-making process disregards the impact of future referrals and assumes that the required number of visits for each referral is equal to her/his expected number of visits. Based on these two features, we developed the Myopic policy to represent current practice. For each new referral, the Myopic policy chooses the best choice (among rejection or assignment to a day within the wait-time target) with a minimum estimated cost. The estimated cost for the assignment of a new referral to a given day is equal to the additional cost incurred by scheduling the expected number of visits of that referral starting from that given day.

The SB policy is an adaptation of the heuristic presented in \citep{Demirbilek2019a} to our problem. It simulates future referrals for a given number of days over multiple scenarios. For each scenario, it initializes a set $\cal Y$ containing the new referrals and simulated future referrals. It chooses a referral in $\cal Y$ whose best choice has the minimum estimated cost among all referrals in  $\cal Y$, performs the best choice and removes that referral from $\cal Y$. The best choice (among rejection or assignment to a day within the wait-time target) for each referral is determined using a similar procedure to the Myopic policy. This process is repeated until $\cal Y$ is empty. After doing this for all scenarios, a new referral is rejected if the number of times it has been accepted in all scenarios is less than a threshold value; otherwise, it is assigned to the most frequent day to which it has been assigned across all scenarios. This threshold value is an important parameter for the SB policy that takes considerable time (12 hours to 2 days for each instance considered in our numerical experiments) to tune. More details about these benchmark policies are given in Appendix A. 

In order to evaluate the performance of any policy $\mu$, we estimate the value functions of a number of states under this policy (i.e., $v_\mu(s)$) by simulating the evolution of the system under the policy for a set period of time. The estimated value function for each state is equal to the sum of the discounted costs over the length of the simulation run. There are three key parameters in the simulation: (i) the number of random states (${{I}^{s}}$), (ii) the number of days in the warm-up period used to determine a random state (${{I}^{w}}$), and (iii) the number of days simulated to estimate the value function for a state (${{I}^{d}}$). We let $({{I}^{s}},{{I}^{w}},{{I}^{d}})=\left( 25,20,365 \right)$, and use the Myopic policy for the duration of the warm-up period. In all experiments, we check the significance of the difference between the performances of two competing policies at the 95\% confidence level.

Table \ref{table-simulation-one} summarizes the results of all the experiments for instances with one service type. The first column provides the name of the experiment and the corresponding parameter that is being varied. All the other parameters remain equal to their default values reported in Table \ref{table-general-settings}. The second column defines the value of the specific parameter with the default value highlighted in bold. Columns \enquote{Rejection hours} and \enquote{Diversion hours} provide the daily average time (in hours) associated with rejected referrals and diverted visits, respectively, for each policy. For example, if we reject a referral with an expected number of visits equal to 8 and service time equal to 0.5 hours, then the corresponding rejection hours is equal to 4. The columns \enquote{Travel time (hours)} and \enquote{Tour length (hours)} show the daily average travel time and the daily average tour length for the nurse for each policy. We consider the Myopic policy as a reference policy to calculate the percentage gap for each of the ALP and SB policies. The percentage gap between the estimated value function of a random initial state $s^0$ for a given policy $\mu$ and the Myopic policy is equal to $100\times({v_{Myopic}^{sim}(s^0)-v_{\mu}^{sim}(s^0)})/{v_{Myopic}^{sim}(s^0)}$. The last two columns, i.e., \enquote{Gap\% (SD)} with sub-labels \enquote{ALP} and \enquote{SB}, report the average percentage gap together with its standard deviation for each of the ALP and SB policies. When the average percentage gap for only one policy (among the ALP and SB) is bolded, that policy significantly outperforms the other, and when both values are bold, no significant difference has been found between the two. The Myopic policy is not the best policy in any of the instances considered in this section.

\begin{table}[h]
	\caption{Simulation results for the experiments with one service type}
	\centering
	\footnotesize
	\resizebox{\textwidth}{!}{%
		\fontsize{9}{13}\selectfont
		\begin{tabular}{M{3cm}M{1cm}M{1cm}M{1cm}M{1cm}M{1cm}M{1cm}M{1cm}M{1cm}M{1cm}M{1cm}M{1cm}M{1cm}M{1cm}M{2cm}M{2cm}}
			\hline
			\multirow{2}{*}{Experiment} &
			\multicolumn{1}{l}{\multirow{2}{*}{Levels}} &
			\multicolumn{3}{c}{Rejection hours} &
			\multicolumn{3}{c}{Diversion hours} &
			\multicolumn{3}{c}{Travel time (hours)} &
			\multicolumn{3}{c}{Tour length (hours)} &
			\multicolumn{2}{c}{Gap\% (SD)} \\ \cline{3-16}
			& \multicolumn{1}{l}{} & ALP  & Myopic & SB & ALP  & Myopic & SB & ALP  & Myopic & SB & ALP  & Myopic & SB & ALP & SB \\ \hline
			\multirow{3}{*}{\parbox{1\linewidth}{\centering $\bar{D}$ (hours) with fixed referral rates}}  & 8                     & 2.04 & 0.18 & 1.59 & 0.06 & 1.56 & 0.36 & 1.33 & 1.47 & 1.51 & 7.29 & 7.73 & 7.48 & \textbf{27.22 (0.93)} & 25.19 (0.97)
			\\  
			& \textbf{8.5} &        2.40 & 0.29 & 1.97 & 0.08 & 1.87 & 0.40 & 1.31 & 1.46 & 1.50 & 7.39 & 7.74 & 7.54 & \textbf{28.60 (0.96)} & 25.62 (1.10)
			\\  
			& 9  &                   2.76 & 0.40 & 2.43 & 0.12 & 2.26 & 0.40 & 1.28 & 1.43 & 1.46 & 7.46 & 7.77 & 7.53 & \textbf{31.72 (0.81)} & 28.43 (0.94)
			\\ \hline
			\multirow{3}{*}{\parbox{1\linewidth}{\centering $\bar{D}$ (hours) with random referral rates}} & 8                     & 2.07 & 0.19 & 1.58 & 0.07 & 1.57 & 0.39 & 1.32 & 1.43 & 1.50 & 7.32 & 7.72 & 7.51 & \textbf{25.93 (1.23)} & 23.81 (1.13)
			\\ 
			& 8.5                   & 2.37 & 0.29 & 1.90 & 0.10 & 1.82 & 0.43 & 1.29 & 1.40 & 1.48 & 7.39 & 7.75 & 7.55 & \textbf{26.81 (0.82)} & 24.91 (0.93)
			\\ 
			& 9                     & 2.72 & 0.39 & 2.37 & 0.15 & 2.25 & 0.46 & 1.26 & 1.38 & 1.44 & 7.48 & 7.77 & 7.57 & \textbf{30.40 (0.79)} & 25.97 (0.87)
			\\ \hline
			\multirow{3}{*}{\parbox{1\linewidth}{\centering Deterministic number of visits}}        & 6                     & 2.23 & 0.27 & 2.19 & 0.12 & 1.86 & 0.03 & 1.28 & 1.45 & 1.39 & 7.49 & 7.79 & 7.66 & 31.68 (0.55) & \textbf{39.11 (0.32)}
			\\ 
			& 8                     & 2.43 & 0.30 & 2.30 & 0.09 & 1.89 & 0.04 & 1.30 & 1.45 & 1.41 & 7.40 & 7.79 & 7.62 & 27.96 (0.78) & \textbf{36.04 (0.46)}
			\\ 
			& 10                    & 2.22 & 0.29 & 2.24 & 0.13 & 1.79 & 0.03 & 1.30 & 1.46 & 1.46 & 7.45 & 7.79 & 7.60 & 28.87 (0.71) & \textbf{34.57 (0.59)}
			\\ \hline
			\multirow{3}{*}{\parbox{1\linewidth}{\centering Poisson number of visits}}             & 6                     & 2.16 & 0.24 & 1.84 & 0.13 & 1.89 & 0.45 & 1.29 & 1.45 & 1.49 & 7.49 & 7.75 & 7.57 & \textbf{33.34 (0.82)} & 26.66 (0.69)
			\\ 
			& \textbf{8}            & 2.40 & 0.29 & 1.97 & 0.08 & 1.87 & 0.40 & 1.31 & 1.46 & 1.50 & 7.39 & 7.74 & 7.54 & \textbf{28.60 (0.96)} & 25.62 (1.10)
			\\ 
			& 10                    & 2.25 & 0.29 & 1.98 & 0.13 & 1.83 & 0.37 & 1.30 & 1.46 & 1.50 & 7.44 & 7.75 & 7.55 & \textbf{29.59 (1.43)} & 25.85 (1.31)
			\\ \hline
			\multirow{3}{*}{\parbox{1\linewidth}{\centering Uniform number of visits}}             & 6 & 2.19 & 0.23 & 1.64	 & 0.12 & 1.91 & 0.66 & 1.30 & 1.46 & 1.53 & 7.49 & 7.75 & 7.59 & \textbf{33.84 (0.90)} & 22.50 (0.98)                    
			\\ 
			& 8 & 2.41	 & 0.25 & 2.28 & 0.08 & 1.88 & 0.32 & 1.32 & 1.46 & 1.57 & 7.38 & 7.76 & 7.28	 & \textbf{29.00 (0.98)}	& 20.00 (1.35)        
			\\ 
			& 10 & 2.29 &	0.28 &	1.66	& 0.18 &	1.98 &	0.75	& 1.29 & 1.45 & 1.51	 & 7.44 & 7.77 & 7.59 & \textbf{29.18 (0.93)}	& 20.46 (0.87)                   
			\\ \hline
			\multirow{4}{*}{\parbox{1\linewidth}{\centering Service time (hours)}}                  & 0.25                  & 2.42 & 0.13 & 2.02 & 0.05 & 2.18 & 0.34 & 1.60 & 1.70 & 1.69 & 7.76 & 7.91 & 7.78 & \textbf{37.43 (0.55)} & 33.58 (0.63)
			\\ 
			& \textbf{0.5}          & 2.40 & 0.29 & 1.97 & 0.08 & 1.87 & 0.40 & 1.31 & 1.46 & 1.50 & 7.39 & 7.74 & 7.54 & \textbf{28.60 (0.96)} & 25.62 (1.10)          \\ 
			& 0.75                  & 2.38 & 0.43 & 2.06 & 0.12 & 1.66 & 0.40 & 1.13 & 1.29 & 1.35 & 7.12 & 7.57 & 7.22 & \textbf{23.85 (1.21)} & 19.41 (1.32)
			\\ 
			& 1                     & 2.50 & 0.64 & 2.36 & 0.18 & 1.70 & 0.41 & 0.98 & 1.15 & 1.22 & 6.86 & 7.30 & 6.86 & \textbf{21.63 (1.60)} & 15.47 (1.39)
			\\ \hline
			\multirow{4}{*}{\parbox{1\linewidth}{\centering Care pattern (days)}}                   & \textbf{1}            & 2.40 & 0.29 & 1.97 & 0.08 & 1.87 & 0.40 & 1.31 & 1.46 & 1.50 & 7.39 & 7.74 & 7.54 & \textbf{28.60 (0.96)} & 25.62 (1.10)          \\ 
			& 3                     & 2.46 & 0.03 & 1.85 & 0.08 & 2.03 & 0.27 & 1.31 & 1.41 & 1.31 & 7.40 & 7.79 & 7.65 & 24.96 (1.26) & \textbf{34.54 (0.86)}
			\\ 
			& 5                     & 2.38 & 0.00 & 2.01 & 0.09 & 1.99 & 0.04 & 1.30 & 1.41 & 1.28 & 7.33 & 7.76 & 7.61 & 19.93 (1.40) & \textbf{37.00 (0.77)}
			\\ 
			& 7                     & 2.25 & 0.00 & 1.95 & 0.12 & 1.86 & 0.05 & 1.30 & 1.40 & 1.31 & 7.30 & 7.66 & 7.55 & 10.64 (1.50) & \textbf{29.70 (1.16)}
			\\ \hline
			\multirow{3}{*}{\parbox{1\linewidth}{\centering Wait-time target (days)}}               & 1                     & 2.66 & 0.3 & 2.77 & 0.05 & 2.08 & 0.03 & 1.37 & 1.47 & 2.14 & 7.24 & 7.66 & 7.19 & \textbf{35.92 (0.67)} & 33.59 (0.67)
			\\ 
			& 3                     & 2.45 & 0.27 & 2.64 & 0.10 & 2.01 & 0.14 & 1.3 & 1.45 & 1.66 & 7.40 & 7.77 & 7.18 & \textbf{32.59 (1.12)} & 26.71 (1.23)
			\\ 
			& \textbf{5}            & 2.40 & 0.29 & 1.97 & 0.08 & 1.87 & 0.40 & 1.31 & 1.46 & 1.50 & 7.39 & 7.74 & 7.54 & \textbf{28.60 (0.96)} & 25.62 (1.10)         \\ \hline
			\multirow{4}{*}{\parbox{1\linewidth}{\centering Diversion weight}}                      & 7.5                   & 2.29 & 0.05 & 1.88 & 0.13 & 2.11 & 0.47 & 1.29 & 1.44 & 1.49 & 7.44 & 7.76 & 7.59 & \textbf{12.71 (1.47)} & \textbf{13.36 (1.05)}
			\\ 
			& \textbf{10}           & 2.40 & 0.29 & 1.97 & 0.08 & 1.87 & 0.40 & 1.31 & 1.46 & 1.50 & 7.39 & 7.74 & 7.54 & \textbf{28.60 (0.96)} & 25.62 (1.10)          \\ 
			& 20                    & 2.39 & 1.25 & 2.65 & 0.07 & 0.96 & 0.13 & 1.32 & 1.47 & 1.56 & 7.38 & 7.70 & 6.96 & \textbf{42.59 (1.25)} & 32.41 (1.55)
			\\ 
			& 1000                  & 2.76 & 1.90 & 2.79 & 0.02 & 0.50 & 0.12 & 1.38 & 1.46 & 1.56 & 7.24 & 7.63 & 6.93 & \textbf{89.11 (0.91)} & 68.18 (1.81)
			\\ \hline
			\multirow{2}{*}{\parbox{1\linewidth}{\centering Overtime limit (hours)}}    & \textbf{0}            & 2.40 & 0.29 & 1.97 & 0.08 & 1.87 & 0.40 & 1.31 & 1.46 & 1.50 & 7.39 & 7.74 & 7.54 & \textbf{28.60 (0.96)} & 25.62 (1.10)          \\ 
			& 2                     & 0.00 & 0.07 & 1.07 & 0.50 & 0.75 & 0.13 & 1.52 & 1.60 & 1.60 & 9.41 & 9.11 & 8.73 & \textbf{21.84 (1.93)} & 12.96 (1.68)
			\\ \hline
			\multirow{3}{*}{\parbox{1\linewidth}{\centering Area diameter (hours)}}                 & \textbf{0.5}          & 2.40 & 0.29 & 1.97 & 0.08 & 1.87 & 0.40 & 1.31 & 1.46 & 1.50 & 7.39 & 7.74 & 7.54 & \textbf{28.60 (0.96)} & 25.62 (1.10)          \\ 
			& 1                     & 3.05 & 0.99 & 3.12 & 0.23 & 2.29 & 0.4 & 2.24 & 2.56 & 2.53 & 7.59 & 7.83 & 7.57 & \textbf{32.25 (0.46)} & 25.44 (0.53)
			\\ 
			& 2                     & 3.88 & 2.81 & 4.71 & 0.49 & 1.97 & 0.19 & 3.38 & 3.97 & 3.63 & 7.59 & 7.74 & 7.18 & \textbf{25.91 (0.49)} & 21.05 (0.56)
			\\ \hline
		\end{tabular}%
	}\label{table-simulation-one}
\end{table}

In the \enquote{$\bar{D}$ (hours) with fixed referral rates} experiment, it is assumed that the referral rates for all regions are the same and equal to $\bar{\lambda }$, and that the varying parameter is the average daily new demand, $\bar{D}$, in hours. Three $\bar{D}$ values of 8.0, 8.5, and 9.0 hours are evaluated. In this experiment, the $\bar{D}$ changes by varying the referral rate $\bar{\lambda }$. The ALP policy outperforms both the Myopic and the SB policies. Almost the same results can be observed for the \enquote{$\bar{D}$ (hours) with random referral rates} experiment, in which the referral rates in different regions are randomly generated.

In the \enquote{Deterministic number of visits} experiment, the number of visits needed for each new referral is assumed to be known, and in the \enquote{Poisson number of visits} and \enquote{Uniform number of visits} experiments, it is assumed to follow a truncated Poisson distribution and a discrete uniform distribution, respectively. Three instances with 6, 8, and 10 visits per new referral (on average for the latter two experiments) are considered in each of these experiments. While the SB policy outperforms the ALP policy in the experiment with a known number of visits per referral, the ALP policy performs best in the more realistic setting with an uncertain number of visits per referral. Furthermore, in the two scenarios involving an uncertain number of visits per referral, the gap between the ALP policy and the SB policy is larger, when the number of visits follows a discrete uniform distribution whose variance is greater than that of the alternative distribution. The SB policy initially schedules the expected number of visits for each referral and adapts as needed. This approach causes this policy to be less effective than the ALP policy in managing undesirable diversions, ultimately resulting in higher overall costs, especially when the number of visits for a referral exceeds the initially anticipated quantity. Nevertheless, if the precise number of visits is known, this strategy offers an advantage over the ALP policy, which only schedules the next immediate visit.

The \enquote{Service time (hours)} experiment evaluates the effect of changing the service time per visit. Four different values of 0.25, 0.50, 0.75, and 1.00 hour are considered. The ALP policy performs better than both the Myopic and the SB policies in these instances, and the average percentage gap between the ALP and Myopic policies increases as the service time decreases.

The \enquote{Care pattern (days)} experiment examines the performance for problem settings with different numbers of days between consecutive visits to a patient, $h_1$. The results for four different values of 1, 3, 5, and 7 show that the average percentage gap between the ALP and Myopic policies decreases with an increase in $h_1$. The ALP policy outperforms the SB policy when $h_1=1$ but is dominated by this policy for the larger $h_1$ values. We will see further in this section (Table \ref{table-K=2}) that for more complex instances with multiple service types and allowable overtime, the ALP policy performs better than the SB policy in most instances with care pattern values larger than one. We should mention that due to discounting, the larger $h_1$, the smaller the difference between the rejection cost of a new referral and the diversion cost of all visits associated with that referral. This would mitigate the consequences of making unfavorable decisions to prevent diversions, particularly concerning the Myopic and SB policies. Furthermore, as $h_1$ increases, the planning horizon for both the Myopic and SB policies becomes significantly longer than that of the ALP policy (e.g., for the instance with $h_1=7$, the lengths of the planning horizon for the ALP, Myopic and SB policies are 7, 54 and 61 days, respectively). However, this significantly increases the execution times of the Myopic and SB policies (e.g., approximately 16 and 8 times longer, respectively, than that of the ALP policy for the instance with $h_1=7$) but it assists them in evaluating various decisions by explicitly observing the implications of scheduling all visits up to the expected number for each referral in a more thorough and extended look-ahead manner. It is important to note that incorporating additional periods to approximate future costs in the ALP approach could lead to tractability issues. 

In the \enquote{Wait-time target (days)} experiment, we evaluate wait-time targets of 1, 3, and 5 days. The ALP policy outperforms both of the other policies. In the \enquote{Diversion weight} experiment, we investigate the effect of changing the weight of the diversion cost, ${{\zeta }^{z}}$, in the cost function \eqref{cost}, while keeping ${{\zeta }^{r}}$, ${{\zeta }^{u}}$ and ${{\zeta }^{q}}$ equal to their default values. The larger ${{\zeta }^{z}}$, the larger the average percentage gap for both the ALP and SB policies. Moreover, the difference between the average percentage gaps for the ALP and SB policies increases in favor of the ALP policy with increasing ${{\zeta }^{z}}$. Some companies may prefer to avoid diversions and therefore would set the diversion cost very high.

In the \enquote{Overtime limit (hours)} experiment, we consider settings with no overtime and with an overtime limit of 2 hours. When overtime is allowed, the average daily overtime for the ALP, Myopic, and SB policies are 1.46, 1.19, and 0.84 hours, respectively. The average percentage gap between the ALP and Myopic policies decreases with the introduction of overtime. This is intuitively reasonable as the main reason for the superiority of the ALP over the Myopic is a reduction in diversions, and the availability of overtime reduces the need for diversions. On the other hand, the ALP policy takes greater advantage of the available overtime capacity to reduce rejections than the SB policy. For the \enquote{Area diameter (hours)} experiment, the ALP outperforms the Myopic and SB policies in all instances with area diameters of 0.5, 1.0, and 2.0 hours. 

Table \ref{table-K=2} presents the results of 25 random instances for the most general case of our problem with two service types, an uncertain number of visits per referral, and overtime. In these instances, $h_k$, $\EX [\hat J_k]$, $e_k$, and $T_k,\forall k$ have been randomly generated based on the set of values defined in Table \ref{table-general-settings}. We report these input parameters in Section \ref{sec:E_add_num}. Other input parameters are at their default values reported in Table \ref{table-general-settings} for $K>1$. The ALP policy is the best choice in 92\% of instances, while the SB policy is the best in only 36\% of them. Thus, the ALP policy is more successful than the other two policies in balancing the rejections, diversions, and overtime to minimize the total cost in the most realistic case of our problem. 

\begin{table}[h]    
	\caption{Simulation results for the most general instances with two service types, an uncertain number of visits per referral, and overtime}
	\centering
	\footnotesize
	\resizebox{\textwidth}{!}{%
		\fontsize{9}{13}\selectfont
		\begin{tabular}{M{1cm}M{1cm}M{1cm}M{1cm}M{1cm}M{1cm}M{1cm}M{1cm}M{1cm}M{1cm}M{1cm}M{1cm}M{1cm}M{1cm}M{1cm}M{1cm}M{2cm}M{2cm}}
			\hline
			\multirow{2}{*}{Ins \#} &
			\multicolumn{3}{c}{Rejection hours} &
			\multicolumn{3}{c}{Diversion hours} &
			\multicolumn{3}{c}{Overtime (hours)} &
			\multicolumn{3}{c}{Travel time (hours)} &
			\multicolumn{3}{c}{Tour length (hours)} &
			\multicolumn{2}{c}{Gap\% (SD)} \\ \cline{2-18}
			& ALP  & Myopic & SB & ALP  & Myopic & SB & ALP  & Myopic& SB & ALP  & Myopic& SB & ALP  & Myopic& SB & ALP & SB              \\ \hline
			1 & 0.96 & 0.19 & 1.33 & 0.61 & 1.25 & 0.46 & 1.26 & 1.43 & 1.14 & 1.46 & 1.60 & 2.21 & 9.07 & 9.32 & 9.00 & \textbf{13.15 (1.07)} & 9.81 (1.22) \\
			2 & 1.12 & 0.01 & 1.27 & 0.28 & 1.26 & 0.41 & 1.57 & 1.73 & 1.46 & 1.52 & 1.62 & 1.78 & 9.52 & 9.71 & 9.44 & \textbf{23.52 (1.27)} & 9.08 (1.49) \\
			3 & 1.35 & 0.03 & 1.50 & 0.39 & 1.32 & 0.36 & 1.16 & 1.41 & 1.01 & 1.27 & 1.38 & 1.63 & 8.88 & 9.35 & 8.88 & \textbf{12.71 (2.11)} & \textbf{11.41 (2.19)} \\
			4 & 0.78 & 0.16 & 1.79 & 0.67 & 1.36 & 0.26 & 1.46 & 1.47 & 1.07 & 1.32 & 1.49 & 1.95 & 9.34 & 9.40 & 8.92 & \textbf{21.39 (1.70)} & 16.24 (1.49) \\
			5 & 1.08 & 0.07 & 1.60 & 0.76 & 1.59 & 0.5 & 1.26 & 1.42 & 1.14 & 1.38 & 1.52 & 2.47 & 8.99 & 9.29 & 8.93 & \textbf{15.87 (1.16)} & 10.83 (1.26) \\
			6 & 1.18 & 0.79 & 1.74 & 0.42 & 0.80 & 0.27 & 1.19 & 1.21 & 0.84 & 1.20 & 1.25 & 1.45 & 8.97 & 9.05 & 8.63 & \textbf{11.32 (1.38)} & 5.32 (1.72) \\	
			7 & 1.01 & 0.00 & 1.41 & 0.62 & 1.35 & 0.40 & 1.35 & 1.60 & 1.31 & 1.31 & 1.46 & 1.57 & 9.20 & 9.59 & 9.28 & \textbf{8.06 (1.78)} & \textbf{8.88 (1.83)} \\
			8 & 0.94 & 0.08 & 1.22 & 0.61 & 1.57 & 0.70 & 1.40 & 1.38 & 1.18 & 1.40 & 1.55 & 2.33 & 9.25 & 9.27 & 9.03 & \textbf{23.54 (0.99)} & 13.47 (1.02) \\
			9 & 0.00 & 0.01 & 0.09 & 1.26 & 1.27 & 1.18 & 1.43 & 1.46 & 1.43 & 1.24 & 1.36 & 1.43 & 9.25 & 9.39 & 9.36 & -1.06 (2.06) & \textbf{4.90 (1.98)} \\
			10 & 0.39 & 0.01 & 0.19 & 0.88 & 1.32 & 1.18 & 1.36 & 1.38 & 1.35 & 1.19 & 1.33 & 1.77 & 9.12 & 9.23 & 9.19 & \textbf{8.66 (2.85)} & 0.05 (2.52) \\
			11 & 0.65 & 0.03 & 0.92 & 0.65 & 1.36 & 0.74 & 1.57 & 1.57 & 1.40 & 1.40 & 1.55 & 1.59 & 9.51 & 9.56 & 9.38 & \textbf{20.54 (1.65)} & 8.30 (1.53) \\
			12 & 2.30 & 0.56 & 2.75 & 0.08 & 1.43 & 0.12 & 1.04 & 1.25 & 0.64 & 1.39 & 1.42 & 2.09 & 8.73 & 9.09 & 8.20 & \textbf{22.89 (1.26)} & 12.05 (1.42) \\
			13 & 1.09 & 0.05 & 0.20 & 0.51 & 1.29 & 1.17 & 1.07 & 1.24 & 1.22 & 1.16 & 1.28 & 1.49 & 8.81 & 9.14 & 9.12 & \textbf{6.35 (2.81)} & \textbf{5.43 (2.34)} \\
			14 & 2.23 & 0.21 & 0.71 & 0.07 & 1.48 & 1.10 & 0.89 & 1.19 & 1.12 & 1.31 & 1.38 & 2.23 & 8.33 & 8.86 & 8.76 & \textbf{14.09 (1.47)} & 2.41 (1.50) \\
			15 & 1.40 & 0.02 & 0.91 & 0.42 & 1.42 & 0.81 & 1.14 & 1.38 & 1.20 & 1.24 & 1.36 & 1.61 & 8.84 & 9.30 & 9.08 & \textbf{13.05 (2.27)} & 6.77 (1.84) \\
			16 & 2.27 & 0.61 & 2.25 & 0.06 & 1.31 & 0.23 & 0.94 & 1.18 & 0.80 & 1.31 & 1.35 & 1.95 & 8.52 & 8.96 & 8.41 & \textbf{22.27 (1.24)} & 15.14 (1.18) \\
			17 & 2.05 & 0.37 & 1.98 & 0.07 & 1.29 & 0.19 & 1.25 & 1.61 & 1.17 & 1.71 & 1.72 & 2.28 & 9.15 & 9.59 & 9.10 & \textbf{22.12 (0.78)} & 17.99 (1.19) \\
			18 & 0.00 & 0.03 & 1.26 & 1.23 & 1.34 & 0.46 & 1.55 & 1.59 & 1.38 & 1.53 & 1.72 & 2.37 & 9.42 & 9.48 & 9.26 & \textbf{7.97 (1.13)} & \textbf{7.54 (1.23)} \\
			19 & 0.79 & 0.00 & 1.32 & 0.53 & 1.21 & 0.32 & 1.63 & 1.76 & 1.42 & 1.49 & 1.63 & 1.67 & 9.59 & 9.75 & 9.40 & \textbf{13.78 (1.71)} & 10.69 (1.89) \\
			20 & 1.94 & 0.30 & 1.87 & 0.06 & 1.26 & 0.20 & 1.31 & 1.64 & 1.25 & 1.70 & 1.69 & 2.22 & 9.24 & 9.62 & 9.20 & \textbf{21.77 (1.25)} & 17.66 (0.98) \\
			21 & 1.70 & 0.51 & 1.52 & 0.19 & 0.98 & 0.25 & 1.11 & 1.35 & 1.09 & 1.35 & 1.40 & 2.01 & 8.87 & 9.26 & 8.97 & \textbf{13.09 (1.11)} & \textbf{14.44 (1.24)} \\
			22 & 0.00 & 0.03 & 0.58 & 1.34 & 1.43 & 1.03 & 1.45 & 1.44 & 1.37 & 1.25 & 1.43 & 1.66 & 9.31 & 9.36 & 9.28 & \textbf{5.07 (2.19)} & \textbf{6.91 (1.67)} \\
			23 & 2.06 & 0.38 & 1.70 & 0.10 & 1.42 & 0.51 & 1.09 & 1.27 & 1.01 & 1.41 & 1.44 & 2.48 & 8.76 & 9.11 & 8.75 & \textbf{22.84 (1.02)} & 11.20 (1.20) \\
			24 & 0.00 & 0.10 & 1.63 & 1.15 & 1.33 & 0.26 & 1.72 & 1.61 & 1.23 & 1.47 & 1.64 & 1.63 & 9.69 & 9.59 & 9.20 & 13.66 (1.46) & \textbf{17.56 (1.49)} \\
			25 & 1.50 & 0.11 & 1.77 & 0.41 & 1.77 & 0.56 & 1.30 & 1.29 & 1.03 & 1.43 & 1.55 & 2.48 & 9.08 & 9.13 & 8.79 & \textbf{26.70 (0.94)} & 13.65 (1.25)
			\\ \hline
		\end{tabular}\label{table-K=2}%
	}
\end{table}

Table \ref{table-k2-special} reports the results for five instances with two service types (with identical service type characteristics as the first five instances in Table \ref{table-K=2}) across three distinct scenarios in which (i) the number of visits per referral is uncertain and overtime in not allowable (\enquote{Poisson number of visits and no overtime}), (ii) the number of visits per referral is deterministic and overtime is allowable (\enquote{Deterministic number of visits and overtime}), and (iii) the number of visits per referral is deterministic and overtime is not allowable (\enquote{Deterministic number of visits and no overtime}). The other input parameters are at their default values reported in Table \ref{table-general-settings} for $K>1$. When the number of visits per referral is uncertain and overtime isn't an option, although the ALP policy still outperforms the SB policy in the instances studied, the SB policy shows greater competitiveness compared to the scenario where overtime is permitted. Moreover, in scenarios with a deterministic number of visits per referral, the ALP policy demonstrates competitive performance compared to the SB policy for instances with two service types, in contrast to its performance in those involving only one service type.

\begin{table}[h]
	\caption{Simulation results for three special cases of instances with two service types}
	\centering
	\footnotesize
	\resizebox{\textwidth}{!}{%
		\fontsize{9}{13}\selectfont
		\begin{tabular}{M{3cm}M{1cm}M{1cm}M{1cm}M{1cm}M{1cm}M{1cm}M{1cm}M{1cm}M{1cm}M{1cm}M{1cm}M{1cm}M{1cm}M{1cm}M{1cm}M{1cm}M{2cm}M{2cm}}
			\hline
			\multirow{2}{*}{Experiment} &
			\multicolumn{1}{l}{\multirow{2}{*}{Ins \#}} &
			\multicolumn{3}{c}{Rejection hours} &
			\multicolumn{3}{c}{Diversion hours} &
			\multicolumn{3}{c}{Overtime (hours)} &
			\multicolumn{3}{c}{Travel time (hours)} &
			\multicolumn{3}{c}{Tour length (hours)} &
			\multicolumn{2}{c}{Gap\% (SD)} \\ \cline{3-19}
			& \multicolumn{1}{l}{} & ALP & Myopic & SB & ALP  & Myopic & SB & ALP  & Myopic & SB & ALP  & Myopic & SB & ALP  & Myopic & SB & ALP & SB \\ \hline
			\multirow{5}{*}{\parbox{1\linewidth}{\centering Poisson number of visits and no overtime}}  
			& 1 & 2.00 & 0.89 & 2.68 & 1.06 & 2.11 & 0.65 & - & - & - & 1.29 & 1.46 & 2.05 & 7.48 & 7.69 & 7.43 & \textbf{17.40 (0.74)} & \textbf{16.66 (0.57)} \\
			& 2 & 2.85 & 0.12 & 3.37 & 0.24 & 2.76 & 0.07 & - & - & - & 1.29 & 1.39 & 1.71 & 7.69 & 7.87 & 7.59 & \textbf{34.01 (0.51)} & 27.36 (0.62) \\
			& 3 & 2.85 & 0.21 & 2.98 & 0.43 & 2.64 & 0.27 & - & - & - & 1.17 & 1.23 & 1.45 & 7.22 & 7.63 & 7.29 & \textbf{26.28 (1.15)} & \textbf{27.92 (1.02)} \\
			& 4 & 3.09 & 0.40 & 3.33 & 0.25 & 2.74 & 0.15 & - & - & - & 1.13 & 1.30 & 1.71 & 7.44 & 7.71 & 7.41 & \textbf{36.70 (0.94)} & \textbf{35.78 (0.90)} \\
			& 5 & 1.51 & 0.38 & 3.53 & 1.39 & 2.53 & 0.16 & - & - & - & 1.19 & 1.36 & 2.10 & 7.53 & 7.67 & 7.17 & \textbf{19.20 (0.74)} & \textbf{18.28 (0.71)} \\
			\hline
			\multirow{5}{*}{\parbox{1\linewidth}{\centering Deterministic number of visits and overtime}}
			& 1 & 0.97 & 0.31 & 1.56 & 0.64 & 1.17 & 0.36 & 1.27 & 1.45 & 1.10 & 1.46 & 1.61 & 2.22 & 9.10 & 9.37 & 8.94 & \textbf{11.75 (0.58)} & \textbf{11.12 (0.54)} \\
			& 2 & 1.13 & 0.04 & 1.32 & 0.33 & 1.23 & 0.38 & 1.52 & 1.73 & 1.42 & 1.51 & 1.62 & 1.76 & 9.47 & 9.72 & 9.40 & \textbf{21.88 (0.91)} & 10.33 (1.01) \\
			& 3 & 1.36 & 0.06 & 1.58 & 0.34 & 1.18 & 0.20 & 1.14 & 1.41 & 0.97 & 1.27 & 1.38 & 1.62 & 8.83 & 9.35 & 8.84 & 9.82 (1.65) & \textbf{14.25 (1.26)} \\
			& 4 & 0.81 & 0.27 & 1.95 & 0.66 & 1.27 & 0.16 & 1.44 & 1.45 & 1.01 & 1.30 & 1.46 & 1.91 & 9.30 & 9.38 & 8.86 & \textbf{19.61 (0.72)} & \textbf{18.10 (0.72)} \\
			& 5 & 1.56 & 0.11 & 1.67 & 0.44 & 1.58 & 0.47 & 1.16 & 1.41 & 1.15 & 1.36 & 1.52 & 2.47 & 8.86 & 9.26 & 8.93 & \textbf{19.20 (0.62)} & 11.54 (0.63) \\ 
			\hline
			\multirow{5}{*}{\parbox{1\linewidth}{\centering Deterministic number of visits and no overtime}}        
			& 1 & 2.20 & 1.23 & 3.52 & 0.97 & 1.83 & 0.08 & - & - & - & 1.31 & 1.46 & 1.87 & 7.45 & 7.69 & 7.28 & 12.68 (0.34) & \textbf{14.16 (0.46)} \\
			& 2 & 2.88 & 0.22 & 3.44 & 0.25 & 2.71 & 0.06 & - & - & - & 1.29 & 1.40 & 1.70 & 7.69 & 7.89 & 7.56 & \textbf{33.42 (0.30)} & 27.20 (0.39) \\
			& 3 & 2.97 & 0.53 & 3.39 & 0.56 & 2.63 & 0.20 & - & - & - & 1.15 & 1.24 & 1.47 & 7.28 & 7.66 & 7.29 & 24.43 (0.46) & \textbf{27.65 (0.42)} \\
			& 4 & 3.10 & 0.50 & 3.41 & 0.23 & 2.66 & 0.08 & - & - & - & 1.13 & 1.29 & 1.68 & 7.44 & 7.72 & 7.41 & \textbf{36.13 (0.40)} & \textbf{35.78 (0.34)} \\
			& 5 & 2.92 & 0.56 & 3.66 & 0.43 & 2.53 & 0.16 & - & - & - & 1.20 & 1.35 & 2.10 & 7.38 & 7.68 & 7.19 & \textbf{26.84 (0.31)} & 19.14 (0.29) \\
			\hline
		\end{tabular}%
	}\label{table-k2-special}
\end{table}

The primary factor contributing to the superiority of the ALP policy compared to the Myopic policy lies in its ability to consider the future consequences of present decisions. While the SB policy also takes this into account by a scenario-based strategy, the effectiveness of this strategy diminishes when handling complex situations. Moreover, the ALP policy's overall average travel time for all instances considered is smaller than those of the Myopic and SB policies by 8.61\% and 22.30\%, respectively. To provide some insights into the tractability of the ALP approach for instances with a higher number of service types, we solved a single instance per $K$, $K \in \{3,4,6,8,10\}$ (because of time restrictions), and reported the results in Section \ref{sec:E_add_num}.

The average execution times needed to tune the ALP (Section \ref{sec:ALP_tuning}) for the instances with 1, 2, 3, 4, 6, 8, and 10 service types are 21.74 minutes, 40.19 minutes, 1.54 hours, 4.61 hours, 15.42 hours, 1.87 days, and 5.28 days, respectively. The equivalent values to derive the ALP parameters using the ALP-1D-2I are 6.60, 30.30, 110.20, 271.84, 785.90, 2661.38, and 4381.25 seconds. The ALP parameters are mainly dependent on factors such as the average referral rates from different regions and the characteristics of service types, and it is reasonable to assume that these parameters would not undergo frequent changes over a short period of time. Once these parameters are derived, we can generate actions on a daily basis using the resulting ALP policy (while those influential factors remain unchanged) in 0.05 to 0.50 seconds for instances with up to 10 service types. The average execution times of a daily action generation using the Myopic and SB policies are significantly larger than that of the ALP policy and range from 0.15 to 3.50 seconds and from 0.12 to 0.87 seconds, respectively. It is worth mentioning that a large portion of the execution time for each of the Myopic and SB policies is dedicated to finding an initial optimal tour for the current visits assigned to each day in the planning horizon. Generally, when the number of rejections is higher, the search for the optimal route becomes focused on a smaller set of locations. Given that the Myopic policy tends to have fewer rejections compared to the SB policy, it consequently requires more execution time.

In Appendix B, we explain a method based on a complete information relaxation by which we can obtain an estimate of an upper bound for the optimality gap of the ALP policy. However, this method is intractable for our default instance even with one service type. Thus, we applied it to a smaller instance with $K=1$, $e_1=1$, $T_1=1$, $\EX [\hat{J}_1]=4$, and $J_1=8$, and the resulting optimality gap was 17.64\%.

\subsection{Illustration of the ALP policy properties}\label{sec:manage}
In this section, we illustrate the consistency of the rejection and assignment actions generated by the ALP policy. This consistency is mainly based on Propositions \ref{Prop_radius} and \ref{Prop_assignment}, and holds for all instances solved in Section \ref{sec:ALPvsMyopic} provided that the ALP parameters are obtained by the ALP-1D-2I (our default method). Consider an instance with a circular area including 60 regions, 2 service types ($h_1=4$, $h_2=2$, $e_1=0.75$, $e_2=0.5$, $\EX [\hat{J}_1]=10$, $\EX [\hat{J}_2]=8$, $T_1=3$, $T_2=5$, $R_1=37.5$, $R_2=20$, $Z_1=7.5$, $Z_2=5$), an area diameter of 2 hours, an average daily new demand of 15 hours ($\bar{D}=15$), no overtime, and the other input parameters at their default values reported in Table \ref{table-general-settings} for $K>1$. The average percentage gap and its standard deviation (in parentheses) for the ALP and SB with respect to the Myopic are respectively 15.79\% (0.39\%) and 13.88\% (0.45\%) for this setting, where the difference between the ALP and each of the Myopic and SB policies is statistically significant. The (\enquote{Rejection hours}, \enquote{Diversion hours}, \enquote{Tour length (h)}) for the ALP, Myopic and SB are respectively (10.21, 0.34, 7.23), (7.12, 3.31, 7.67) and (10.17, 0.50, 7.53).

\begin{figure}[h] 
	\begin{center}
		\includegraphics[height=2.5in]{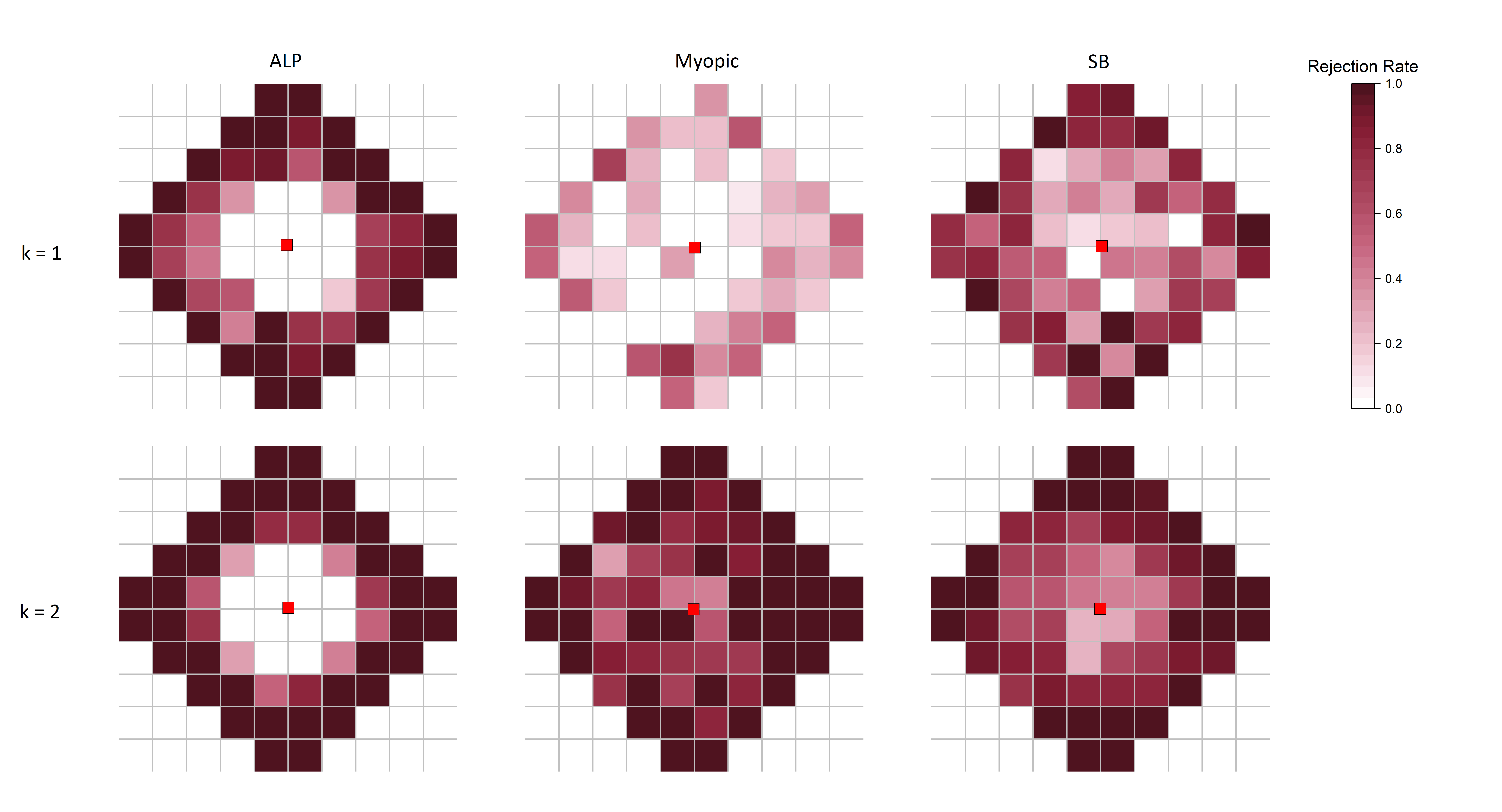}
		\caption{Rejection rates for different policies, regions and service types} \label{fig-Rejection}
	\end{center}
\end{figure}

Figure \ref{fig-Rejection} illustrates the resulting rejection rate for the different policies, regions and service types. The rejection rates are computed as the number of rejections over the total number of referrals for each region and service type. The darker a region, the higher the rejection rate.

For the ALP policy, there are two circles around the depot for each service type. Inside the first circle, called the \enquote{sure acceptance} area, all referrals are accepted. The area between the first and the second circles represents the \enquote{maybe} area, where some referrals may be rejected depending on the state of the system. Beyond the second circle is the \enquote{sure rejection} area where all referrals are rejected. In the multi-nurse setting, this result could be used to determine the appropriate size of the coverage area for each nurse and how to overlap the \enquote{maybe} areas associated with different nurses in order to increase the overall service level. Proposition \ref{Prop_radius} helps us determine the radii of these circles based on the values of the ALP parameters. The radii of the first and the second circles for service type 1 in the ALP policy are 2 and 4 distance units (length of each region), respectively. However, the values for service type 2 are 2 and 3 distance units. The ALP policy tends to reject more referrals of type 2 than type 1 that could be explained by ${{R}_{1}}>{{R}_{2}}$.

Similar to the ALP policy, the Myopic and SB policies reject more referrals of type 2 than type 1. However, the differences in the rejection rates between the two service types are more severe for the Myopic policy. For both Myopic and SB policies, rejection is possible in all regions, even those closest to the depot. On average, the rejection rate of each of the Myopic and SB policies increases with the distance from the depot but does not follow a regular pattern.

\begin{figure}[h]
	\begin{center}
		\includegraphics[height=2.5in]{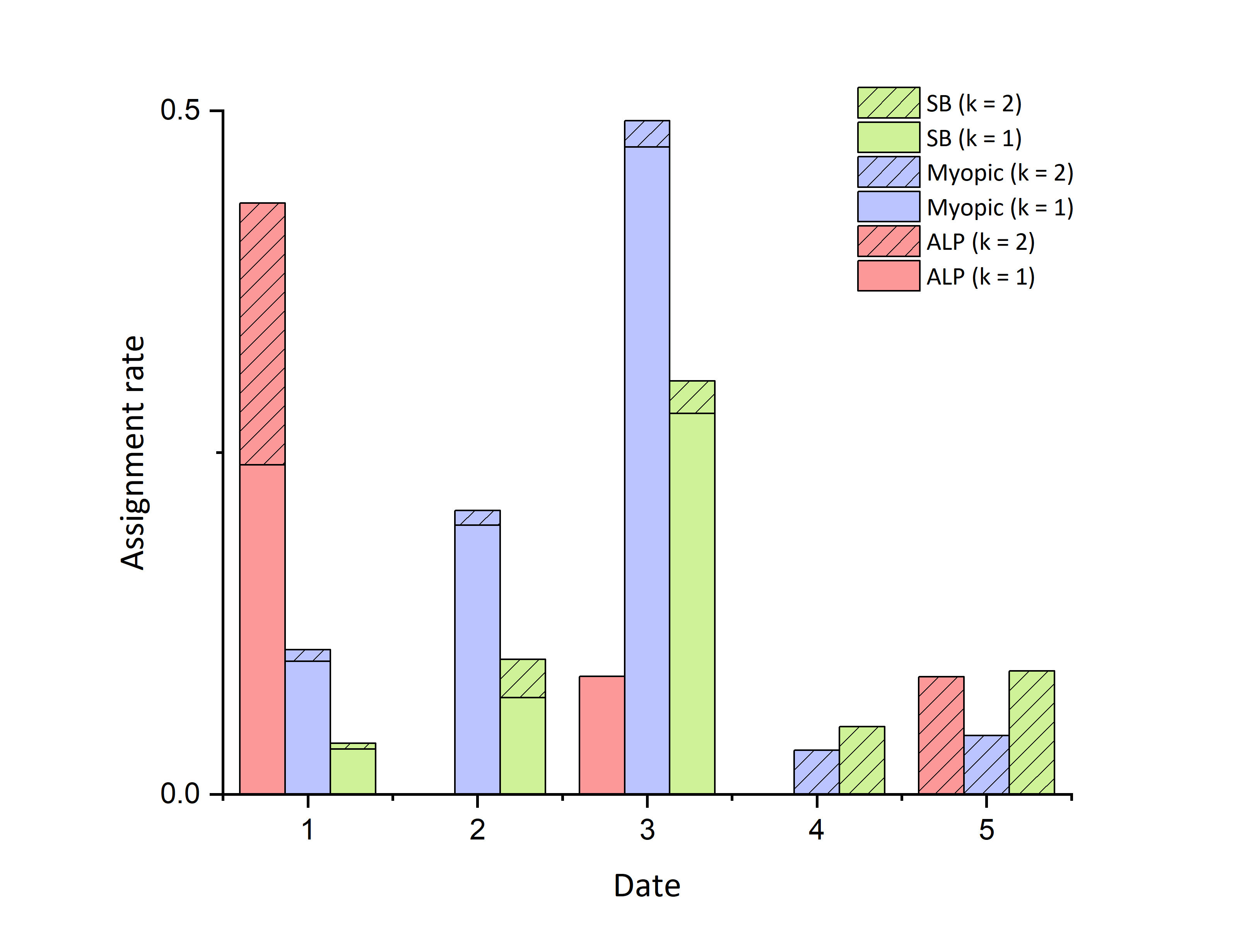}
		\caption{Referral assignment rates for different dates, policies and service types} \label{fig-Assignment}
	\end{center}
\end{figure}

Figure \ref{fig-Assignment} provides the rates at which different days within the wait-time targets are assigned to new referrals for each policy and service type. These rates are computed as the ratio of the number of referrals assigned to a particular date to the total number of referrals of a specific service type. For the ALP policy, each referral is assigned to either the first day or the last day within the wait-time target, as indicated in Proposition \ref{Prop_assignment}. On the other hand, the Myopic and SB policies assign referrals to every day within the wait-time target with higher assignment rates for later days.

\section{Conclusions and suggestions for future research}\label{sec:Conclusions-and-suggestions}
In this study, we focused on the stochastic HRSP. We developed an MDP model for a single-nurse dynamic HRSP with uncertain new referrals and number of visits per referral. We solved this MDP using an ALP approach. Despite the ALP approach, the most general form of the HRSP remained intractable. To overcome this obstacle, we determined a closed-form solution for the optimal ALP parameters of a special case of the problem that in turn helped us derive two heuristic reduction techniques. These techniques allowed our approach to determine the ALP parameters of large-scale instances of the general problem.

We compared the ALP policy with a myopic policy (referred to as Myopic) and a scenario-based policy (referred to as SB). Unlike the Myopic policy, the ALP and SB policies consider the future consequences of present decisions. While the ALP policy schedules only the next immediate visit, the Myopic and SB policies initially schedule the expected number of visits for each referral and adapt as needed. This approach leads to notably longer planning horizons and execution times for these benchmark policies compared to the ALP policy, especially for large care patterns. Moreover, it often causes the benchmark policies to be less effective than the ALP policy in managing undesirable diversions, particularly when the number of visits for a referral exceeds the initially expected quantity. However, it offers an advantage over the ALP policy in some cases, especially when the precise number of visits is known. The ALP policy outperformed the Myopic policy representing current practice, generating better solutions in terms of a discounted weighted sum of rejection, diversion, overtime, and travel time costs. Although the SB policy had a better performance than the ALP policy in a few special cases (e.g., instances with a single service type \textit{and} a deterministic number of visits per referral or a care pattern larger than one), its effectiveness declined in more complex problem settings. In the most general case of our problem, involving multiple service types, overtime, and an uncertain number of visits per referral, the ALP policy outperformed the SB policy. Using a method based on a complete information relaxation, we derive an optimality gap of 17.64\% for a specific instance. The ALP policy strategically rejects demand in order to pre-empt downstream congestion that would lead to a large number of diversions. In addition, it accepts all referrals from patients located inside an inner circle centered at the depot and rejects all referrals from patients located beyond an outer circle centered at the depot, leaving an area in-between where the acceptance/rejection decision depends on the current congestion of the system. This resulted in a significantly smaller travel time for this policy compared with the benchmark policies. Furthermore, the ALP policy always assigns an accepted referral to either the first day or the last day within the corresponding wait-time target. These features make the actions generated by the ALP policy more consistent than those of the benchmark policies. We provided a simple procedure to tune the state-relevance weights in the ALP. While tuning the ALP and deriving the ALP parameters are time-consuming for instances with numerous service types, they only need to be performed once and then actions can be generated very fast on a daily basis using the resulting ALP policy (e.g., in less than one second for instances with 10 service types and significantly faster than the benchmark policies considered).

Some directions for future research include (i) extending the proposed ALP approach to solve a version of the problem involving multiple nurses with different depots (a simple heuristic for this purpose is given in Section \ref{sec:E_multi-nurse}), (ii) including time-windows in the model, (iii) finding faster ways to tune the ALP, (iv) investigating alternative basis functions to be employed in the approximate linear function in order to enhance the ALP policy, and (v) quantifying the fairness of service to all populations in an area and including this metric in the objective to mitigate bias against distant locations.

\section*{Appendices}
\subsection*{Appendix A: More details about the benchmark policies\label{sec:A-myopic-policy}}
In the Myopic policy, if a new referral of type $k$ is accepted, ${{\bar{J}}_{k}}$ visits are initially scheduled for her/him. Therefore, the planning horizon ($T$) for the Myopic policy is equal to $\underset{k}{\mathop{\max }}\,\{{{T}_{k}}+{{h}_{k}}({{\bar{J}}_{k}}-1)\}$. For each existing patient of type $k$, the remaining visits up to her/his ${{\bar{J}}_{k}}$-th visit are scheduled, and if the patient has already been served ${{\bar{J}}_{k}}$ times, only one visit is scheduled for her/him. In order to determine the specific action to take in state $s=(\vec{x},\vec{y})$ under the Myopic policy, we initially find an optimal tour for the current visits assigned to each day in the planning horizon, and then we perform the following algorithm:\vspace{5pt}

\noindent
\textbf{Algorithm 1. Myopic policy}
\begin{steps}
	\item There are $\sum\limits_{kl}{{{y}_{kl}}}$ new referrals. Sort them in non-decreasing order of distance from the depot. Label the $i$-th referral in the resulting list ${{\Pi }_{i}}$, and set $i=1$. The service type for this patient is $k,k\in \{1,...,K\}$.
	\item Make a decision about ${{\Pi }_{i}}$ according to the following procedure. Each new patient can be assigned to a day $t,t\in \{1,...,{{T}_{k}}\}$ or rejected. Calculate the resulting additional cost of each of these ${{T}_{k}}+1$ choices and denote them by ${{\Delta }_{{{t}'}}},{t}'\in \{0,1,...,{{T}_{k}}\}$, where ${{\Delta }_{0}}$ corresponds to the rejection choice and is equal to ${{R}_{k}}={{\zeta }^{r}}\EX[{{\hat{J}}_{k}}]{{e}_{k}}$. To calculate ${{\Delta }_{{{t}'}}}$ corresponding to the assignment of the patient to day ${t}',{t}'\in \{1,...,{{T}_{k}}\}$, let ${t}'=1$ and go to the next step. 
	\item Assign the first visit of ${{\Pi }_{i}}$ to day ${t}'$ and assign the remaining ${{\bar{J}}_{k}}-1$ visits of this patient following her/his care pattern up to at most day $T$. Update the nurse’s route on each day to which a visit of ${{\Pi }_{i}}$ has been assigned. A new visit assigned to each day is scheduled based on the cheapest insertion method \citep{Ulmer2018}. After insertion of the new visit in a day’s route, three scenarios may happen: (a) if the route length is not larger than $\chi $ then the resulting additional cost of assigning the patient’s visit to this day would be equal to the additional travel time cost, (b) if the route length is not larger than $\chi +{\chi }'$ but is larger than $\chi $ then the additional cost to this day would be equal to the sum of the additional travel time cost and the additional overtime cost, or (c) if the route length is larger than $\chi +{\chi }'$ then the additional cost in this day would be equal to the diversion cost of this visit. Calculate the sum of the discounted additional costs for patient ${{\Pi }_{i}}$ and set ${{\Delta }_{{{t}'}}}$ equal to this value. Let ${t}'={t}'+1$. If ${t}'\le {{T}_{k}}$, repeat this step; otherwise, go to the next step.
	\item Select a choice for ${{\Pi }_{i}}$ with the lowest ${{\Delta }_{{{t}'}}},{t}'\in \{0,...,{{T}_{k}}\}$, update the schedule based on this choice, and go to the next step.
	\item Let $i=i+1$. If $i \le \sum\limits_{kl}{{{y}_{kl}}}$, return to Step 2; otherwise, stop. 
\end{steps}

In the SB policy, we have adapted the heuristic presented in \citep{Demirbilek2019a} to our problem. In order to take into account the future referrals in our decision-making, we generate several scenarios for the next $T'$ days by generating random new referrals based on the probability distribution function of the new referrals of each service type in each location. The following algorithm explains the steps to perform the SB policy for a given state: \vspace{5pt}

\noindent
\textbf{Algorithm 2. SB policy}
\begin{steps}
	\item Generate $N_{sc}$ random scenarios for future referrals in the next $T'$ days. Let $sc=1$. 
	\item Initialize a set $\cal Y$ by assigning the new referrals and simulated future referrals in scenario $sc$ to it.
	\item For each referral in $\cal Y$, find the best choice (among rejection or assignment to a day within the wait-time target) according to the procedure given in Steps 2 and 3 of the Myopic algorithm. 
	\item Find a referral in $\cal Y$ whose best choice has the minimum cost among all, perform the best choice for it, and remove it from $\cal Y$. If $\cal Y$ is empty, go to the next step; otherwise back to Step 3.
	\item Let $sc=sc+1$. If $sc \le N_{sc}$, back to Step 2; otherwise proceed to the next step.
	\item Reject each new referral if the number of times it has been accepted in all scenarios is less than a threshold value $N_{tr}$; otherwise, assign the new referral to the most frequent day it has been assigned in all scenarios.
\end{steps}
It is assumed that $T'=\max \{\underset{k}{\mathop{\max }}\,\{{{h}_{k}}\},\underset{k}{\mathop{\max }}\,\{{{T}_{k}}\},5\}$, and the length of the planning horizon, $T$, is equal to $T'+\underset{k}{\mathop{\max }}\,\{{{T}_{k}}+{{h}_{k}}({{\bar{J}}_{k}}-1)\}$. We set $N_{sc}=100$, and we tuned the $N_{tr}$ for each instance by checking different values in the set $\{0,10,20,...,100\}$. The larger the care pattern or the smaller the service time, the longer it takes to tune the $N_{tr}$.

\subsection*{Appendix B: Lower bounds for the optimal value function \label{sec:App-LB}}

In a popular approach to obtaining a lower bound for the optimal value function of a given state, the nonanticipativity constraints are relaxed assuming perfect information about all future stochastic elements (see, e.g., \citealp{Brown20171355}). To fix the issue of simulating infinite sequences in the perfect information relaxation of our discounted infinite horizon problem, a standard and equivalent formulation is used in which there is no discounting, but rather an absorbing, costless state $\bar{s}$ that is reached with probability $1- \gamma$ in each period (see, e.g., \citealp{Puterman1994}, Chapter 5). In each sample of the perfect information relaxation, an absorption time of $\bar{T}$ is revealed (following a geometric distribution with a parameter $1-\gamma$) along with a sample path which provides all information about future referrals up to day $\bar{T}$ including the exact number of visits. The following recursion must be solved for each sample path:
\begin{equation}
	v_t^{\cal G}(s_t)=\underset{a\in A_{s_{t}}}{\mathop{\min }}\,\{c(s_t,a)+ v_{t+1}^{\cal G}(s_{t+1})\},\forall s_t\in S,t\in \{0,...,\bar{T}-1\},\label{Dual_LB_DP}
\end{equation}

where $v_{\bar{T}}^{\cal G}(\bar{s})=0$, and $s_{t+1}$ represents the state $s_t$ after a transition considering the perfect information relaxation. Even though solving recursion \eqref{Dual_LB_DP} is much simpler than the primal DP \eqref{Bel} it is still intractable for our default instance even with one service type (Table \ref{table-general-settings}). Therefore, to have a tractable number of states at different stages of recursion \eqref{Dual_LB_DP}, we considered an instance with one service type in which $\bar D=7.5$ hours, $e_1=1$ hour, $T_1=1$ day, the number of visits for an accepted referral is in the range of $[1-8]$ with a mean value of 4, and other parameters are equal to their default values given in Table \ref{table-general-settings}. Because $T_1=1$, if a referral is accepted, it will definitely be assigned to day 1. The average percentage gap and its standard deviation (in parentheses) of the ALP and SB with respect to the Myopic for this instance were respectively 34.52\% (1.06\%) and 28.56\% (1.18\%). 

Based on the features of this instance and also the complete information in each sample path, the state $s_t$ is represented by: (i) the rejection/acceptance status of each referral with $J$ visits arriving at day $t',t'<t$ for which we have $t'+J-1 \ge t$, (ii) non-zero $x^0_{1,1,lj}$ values in the initial state $s_0=(\vec x^0,\vec y^0)$ for which we have $J-j-1 \ge t$, and (iii) the new referrals waiting for a decision at day $t$. If at most four new referrals come each day (with respect to the average daily new demand considered for the instance), then there can be a maximum of $(2^4)^7\approx2.7e9$ different states at each day. However, the real number of states occurring at each day of each sample path for this instance is much less than this value but it still results in a long computational time. We obtained five random initial states $s_0$ after a warm-up period using the Myopic policy, and generated 100 sample paths for each of them. Then, we solved recursion \eqref{Dual_LB_DP} for each initial state and each of its sample paths. Moreover, we evaluated the ALP policy for each initial state when referrals arrive according to each of its sample paths. The resulting average percentage gap between the lower bound obtained by a complete information relaxation and the upper bound derived by the ALP policy was 17.64\% with a standard deviation of 0.33\%. The average execution time to calculate this lower bound for each initial state was 61.58 hours.

For any ALP parameters set that is feasible in the ALP model \eqref{P-ALP-1}-\eqref{P-ALP-3}, the approximate linear function \eqref{eq:affine_fun} provides a lower bound for the optimal value function of a given state \citep{Desai2013452}.

\begin{proposition}
	The set of the ALP parameters obtained by each of the ALP-1D and ALP-1D-2I is feasible in the ALP model of ${\cal P}_g$. \label{proposition:reduction_feasibility}
\end{proposition}
\noindent\textsc{Proof.} The proof is done for the ALP-1D case but it is easily expandable for the ALP-1D-2I one. Let $\pi({\cal P}^1_g)=(\tilde{\eta }({\cal P}^1_g),{{\tilde{\tau }}_{tklj}}({\cal P}^1_g),{{\tilde{\rho }}_{kl}}({\cal P}^1_g))$ be the ALP parameters set obtained for the proxy ${\cal P}^1_g$, and $\pi({\cal P}_g)$ be the corresponding ALP parameters set for ${\cal P}_g$ using the ALP-1D. Assume that the $\pi({\cal P}_g)$ is infeasible in the ALP model of ${\cal P}_g$. Then, there exists at least one feasible state-action pair $(\hat s,\hat a)$ in ${\cal P}_g$ for which we have $(1-\gamma )\tilde{\eta }({\cal P}^1_g) +\sum\limits_{tk{l^*}j}{{{\tilde{\tau }}_{tk{l^*}j}}({\cal P}^1_g)\sum\limits_{l \in {\cal L}_{l^*}}\left( {{\hat x}_{tklj}}-\gamma \EX[X'_{tklj}(\hat s,\hat a)] \right)}+\sum\limits_{k{l^*}}{{{\tilde{\rho }}_{k{l^*}}}({\cal P}^1_g)\sum\limits_{l \in {\cal L}_{l^*}}\left( {{\hat y}_{kl}}-\gamma \EX[{{Y}_{kl}}]({\cal P}_g) \right)} > c(\hat s,\hat a)$, where $\hat s = (\hat x_{tklj}, \forall tklj,\hat y_{kl}, \forall kl)$. Create an state-action pair $(\hat s^1, \hat a^1)$ for ${\cal P}^1_g$ based on the state-action pair $(\hat s,\hat a)$, where $\hat x^1_{tk{l^*}j}=\sum\limits_{l \in {\cal L}_{l^*}}\hat x_{tklj}, \forall tk{l^*}j$, $\hat y^1_{k{l^*}}=\sum\limits_{l \in {\cal L}_{l^*}}\hat y_{kl}, \forall k{l^*}$, $\hat r^1_{k{l^*}}=\sum\limits_{l \in {\cal L}_{l^*}}\hat r_{kl}, \forall k{l^*}$, $\hat n^1_{tk{l^*}}=\sum\limits_{l \in {\cal L}_{l^*}}\hat n_{tkl}, \forall tk{l^*}$, $\hat z^1_{k{l^*}}=\sum\limits_{l \in {\cal L}_{l^*}}\hat z_{kl}, \forall k{l^*}$, and the optimal order of visiting the patients with a minimum travel time is obtained through typical routing methods. Because the travel time for the state-action pair $(\hat s^1, \hat a^1)$ (in a 1D array) is no greater than that of the state-action pair $(\hat s,\hat a)$ (in a 2D area), $(\hat s^1, \hat a^1)$ is feasible for  ${\cal P}^1_g$. Therefore, we have $(1-\gamma )\tilde{\eta }({\cal P}^1_g) +\sum\limits_{tk{l^*}j}{{{\tilde{\tau }}_{tk{l^*}j}}({\cal P}^1_g)\left( {{\hat x}^1_{tklj}}-\gamma \EX[X'_{tk{l^*}j}(\hat s^1,\hat a^1)] \right)}+\sum\limits_{k{l^*}}{{{\tilde{\rho }}_{k{l^*}}}({\cal P}^1_g)\left( {{\hat y}^1_{kl}}-\gamma \EX[{{Y}_{k{l^*}}}]({\cal P}^1_g) \right)} \le c(\hat s^1, \hat a^1)$ that is in contradiction with our initial assumption because $c(\hat s^1, \hat a^1) \le c(\hat s, \hat a)$. 

The lower bounds obtained by the value of the approximate linear function \eqref{eq:affine_fun} with the ALP parameters derived by the ALP-1D-2I for the same instance and initial states used in the previous method (based on the complete information relaxation) were much smaller than those obtained by this method. Because of that, we have not reported these lower bounds. We have observed a similar performance from this lower bound in \citep{Patrick2008,Klapp2016}. The quality of this lower bound is determined largely by the choice of basis functions in the approximate linear function \citep{Desai2013452,DeFarias2003850}. However, finding a better set of basis functions other than those considered is not straightforward, and even if possible, would be less interpretable and harder to solve.

In another approach to obtaining a lower bound, the nonanticipativity constraints are relaxed, and a penalty is imposed based on an approximate value function to punish violations of the nonanticipativity \citep{Brown20171355}. We applied the approximate linear function to generate penalties in this approach. While the resulting lower bound was better than the lower bound obtained by the approximate linear function, it was worse than the lower bound obtained by the complete information relaxation with zero penalties.

\newcommand{\changeitem}{%
	\let\latexitem\item
	\renewcommand\item[1][]{\latexitem\relax{##1 --} }%
}

\vspace{5pt}
\noindent\textbf{Acknowledgement} \\
This research was partially supported by the Telfer School of Management at University of Ottawa [SMRG Postdoctoral Research Fellowship Support Grant 2020] and by the Natural Sciences and Engineering Research Council of Canada (NSERC) [RGPIN-2018-05225 and RGPIN-2020-210524].

\bibliographystyle{unsrtnat}
\bibliography{references}  






\newpage
\setcounter{page}{1}



\section*{Electronic Companion}
\renewcommand{\thesection}{EC.\arabic{section}}
\setcounter{section}{0}

\renewcommand{\thefigure}{EC.\arabic{figure}}
\setcounter{figure}{0}
\renewcommand{\thetable}{EC.\arabic{table}}
\setcounter{table}{0}

\renewcommand{\thealgorithm}{EC.\arabic{algorithm}}
\section{Tuning the ALP}\label{sec:ALP_tuning}
Unlike the non-approximated LP formulation of an \ac{MDP} (i.e., \eqref{LP}), the ALP does in fact depend on the choice of $\EX_{\alpha}[{{X}_{tklj}}],\forall tklj$ and $\EX_{\alpha}[{{Y}_{kl}}],\forall kl$ values. For simplicity, we assume $\EX_{\alpha}[{{X}_{tklj}}]=\varepsilon ,\forall tklj\notin \Theta $, and zero otherwise, where $\varepsilon $ is an unknown constant. Also, we let ${{E}_{\alpha}}[{{Y}_{kl}}]=\EX[{{Y}_{kl}}],\forall kl$. The optimal \ac{ALP} parameters change as we vary $\varepsilon $. However, the number of different solutions is limited and dependent on the value of specific parameters such as the number of service types and number of regions in the corresponding one-dimensional array for an area (i.e., $K$ and ${{L}^{*}}$). The numerical experiments show a trade-off between the rejection and diversion costs associated with the resulting \ac{ALP} policies for different values of $\varepsilon $. The rejection cost increases and the diversion cost decreases with an increase in $\varepsilon $. The minimum number of rejections happens when $\varepsilon =0$ (though often achieved before $\varepsilon $ reaches zero). The minimum number of rejections is not necessarily zero for systems with large demand. For very large values of $\varepsilon $, all arrivals are rejected, and therefore the minimum number of diversions occurs. 

\begin{figure}[t]
	\begin{center}
		\includegraphics[height=2in]{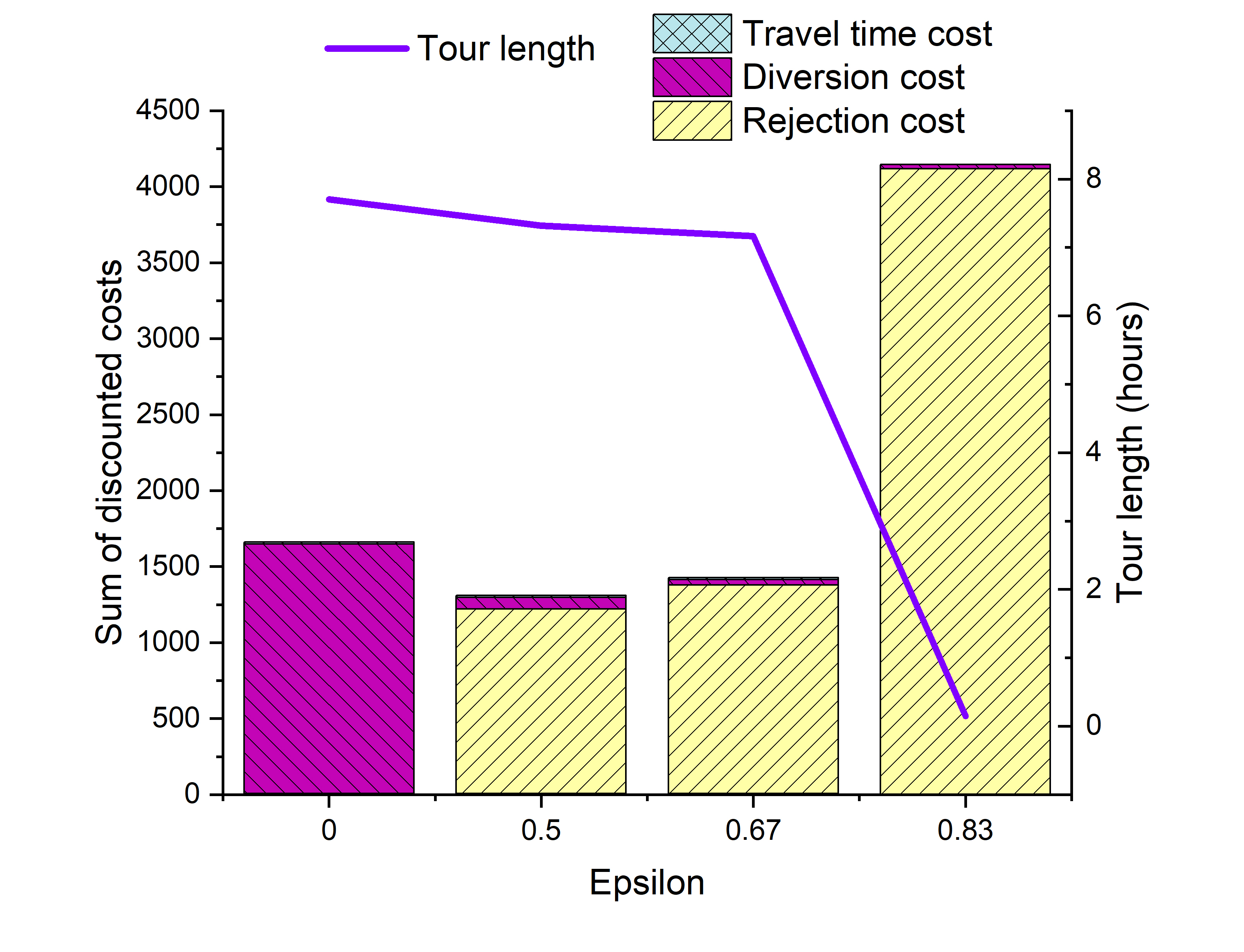}
		\caption{Trade-off between the rejection and the diversion costs associated with the ALP policy as the value of $\varepsilon $ increases} \label{fig-Tradeoff}
	\end{center}
\end{figure}

Consider an example with only one service type and the default parameter values presented earlier in Table \ref{table-general-settings}. Figure \ref{fig-Tradeoff} shows the trade-off between the rejection cost and the diversion cost as the value of $\varepsilon $ increases. The height of each bar for each specific value of $\varepsilon $ represents the estimated average value function for the corresponding \ac{ALP} policy, and the different colors represent the average sum of the discounted rejection, diversion, and travel time costs. When $\varepsilon =0$, there is no rejections and we have a large diversion cost. As $\varepsilon $ increases, more rejections and less diversions are observed. In the most extreme case, when $\varepsilon =0.83$, the resulting \ac{ALP} policy rejects all arrivals. The slight diversion cost seen for this extreme case corresponds to the remaining visits of already assigned patients in the initial state after warm-up. The minimum estimated value function is achieved when $\varepsilon =0.50$.

In the following, an algorithm is presented to systematically adjust $\varepsilon $ for the ALP so as to minimize the average value function estimate. In this algorithm, ${{\pi }_{\varepsilon }}=(\tau ,\rho ,\eta )$ represents the optimal ALP parameters if $\EX_{\alpha}[{{X}_{tklj}}]=\varepsilon ,\forall tklj\notin \Theta $.

\noindent
\textbf{Algorithm EC.1. Tuning $\varepsilon $}
\begin{steps}
	\item Let ${{\varepsilon }_{1}}=0$ and ${{\varepsilon }_{2}}=\sum\nolimits_{kl}{\EX[{{Y}_{kl}}]}/KL$.
	\item While ${{\pi }_{{{\varepsilon }_{1}}}}$ and ${{\pi }_{{{\varepsilon }_{2}}}}$ are the same, let ${{\varepsilon }_{1}}={{\varepsilon }_{2}}$ and ${{\varepsilon }_{2}}=2{{\varepsilon }_{2}}$.
	\item Let ${{\varepsilon }_{cur}}={{\varepsilon }_{1}}$, and calculate the average value function of the ALP policy for ${{\varepsilon }_{cur}}$.
	\item Let ${\varepsilon }'=({{\varepsilon }_{1}}+{{\varepsilon }_{2}})/2$. If ${{\pi }_{{{\varepsilon }'}}}$ and ${{\pi }_{{{\varepsilon }_{2}}}}$ are the same, let ${{\Delta }_{\varepsilon }}={{\varepsilon }_{2}}-{{\varepsilon }_{1}}$, ${{\varepsilon }_{cur}}={{\varepsilon }_{2}}$ and go to Step 5; else if ${{\pi }_{{{\varepsilon }'}}}$ and ${{\pi }_{{{\varepsilon }_{1}}}}$ are the same, then let ${{\varepsilon }_{1}}={\varepsilon }'$ and repeat this step; otherwise, let ${{\varepsilon }_{2}}={\varepsilon }'$, and repeat this step.
	\item If ${{\pi }_{{{\varepsilon }_{cur}}}}$ and ${{\pi }_{{{\varepsilon }_{1}}}}$ are the same, let ${{\varepsilon }_{1}}={{\varepsilon }_{cur}}$, ${{\varepsilon }_{2}}={{\varepsilon }_{cur}}+2{{\Delta }_{\varepsilon }}$, and go back to Step 4; otherwise, calculate the average value function of the ALP policy for ${{\varepsilon }_{cur}}$, and proceed to the next step.
	\item Stop if, under the ALP policy for ${{\varepsilon }_{cur}}$, all arrivals are rejected, and select the ${{\varepsilon }_{cur}}$ resulting in the minimum average value function; otherwise, let ${{\varepsilon }_{1}}={{\varepsilon }_{cur}}$, ${{\varepsilon }_{2}}={{\varepsilon }_{cur}}+{{\Delta }_{\varepsilon }}$, and go back to Step 4.
\end{steps}

In order to speed-up the process, if ${{\pi }_{{{\varepsilon }_{cur}}}}$ and ${{\pi }_{{{\varepsilon }_{1}}}}$ are the same in Step 5, there is no need to calculate the cost of ${{\varepsilon }_{cur}}$ as it is equal to that of ${{\varepsilon }_{1}}$ we have already calculated. The execution times of this algorithm for instances with 1, 2, 3, 4, 6, 8, and 10 service types are respectively 21.74 minutes, 40.19 minutes, 1.54 hours, 4.61 hours, 15.42 hours, 1.87 days, and 5.28 days.

\section{Proof of Theorem \ref{theorem1}}\label{sec:proof-theorem1}
For simplicity in this section, we denote the closed-form solution ${{\eta }^{*}}({\cal P}_s)$, $\tau _{tklj}^{*}({\cal P}_s)$, and $\rho _{kl}^{*}({\cal P}_s)$ obtained through \eqref{T1_1}-\eqref{T1_4} by $\eta $, ${{\tau }_{tklj}}$, and ${{\rho }_{kl}}$, respectively.
The three main steps of the proof are: 
\begin{steps}
	\item Prove the feasibility of the hypothesized closed-form solution for the ALP. 
	\item Determine the necessary and sufficient conditions under which the hypothesized primal solution and a corresponding dual solution would satisfy the complementary slackness conditions.
	\item Demonstrate that there exists a dual solution satisfying the necessary and sufficient conditions obtained in Step 2. 
\end{steps}
In the following, we explain these steps one by one.
\subsection{Step 1}
If the hypothesized solution obtained from \eqref{T1_1}-\eqref{T1_4} is true, then constraint set \eqref{P-ALP-2} must be satisfied for all state-action pairs. Because all referrals are accepted, we have ${{y}_{kl}}=\sum\limits_{t;t\le {{T}_{k}}}{{{n}_{tkl}}}$. Also, from \eqref{T1_4}, we have $(1-\gamma )\eta =\gamma \sum\nolimits_{kl}{{{\rho }_{kl}}\EX[{{Y}_{kl}}]}$. After some manipulations, constraint set \eqref{P-ALP-2} can be written as follows:
\begin{multline}
    \sum\limits_{kl}{{{Z}_{k}}{{z}_{kl}}}+Qq-\sum\limits_{klj}{\left( {{\tau }_{1,klj}}-\gamma {{p}_{k,j+2}}{{\tau }_{{{h}_{k}},kl,j+1}} \right){{x}_{1,klj}}}-\sum\limits_{tklj;2\le t<{{h}_{k}},j\ne 0}{\left( {{\tau }_{tklj}}-\gamma {{\tau }_{t-1,klj}} \right){{x}_{tklj}}}-\\
	\sum\limits_{tklj;2\le t<{{T}_{k}},j=0}{\left( {{\tau }_{tklj}}-\gamma {{\tau }_{t-1,klj}} \right){{x}_{tklj}}}+ \gamma \sum\limits_{kl}{{{p}_{k,2}}{{\tau }_{{{h}_{k}},kl,1}}{{n}_{1,kl}}}+\gamma \sum\limits_{tkl;2\le t\le {{T}_{k}}}{{{\tau }_{t-1,kl,0}}{{n}_{tkl}}}- \\
    \sum\limits_{tkl;t\le {{T}_{k}}}{{{\rho }_{kl}}{{n}_{tkl}}}\ge 0, \forall s\in S,\forall a\in {{A}_{s}}. 
\end{multline}
We know that ${{\tau }_{tklj}}-\gamma {{\tau }_{t-1,klj}}=0,\forall tklj,t\ge 2$. Also, from \eqref{T1_1}, we also have:
\[{{\tau }_{1,klj}}-\gamma {{p}_{k,j+2}}{{\tau }_{{{h}_{k}},kl,j+1}}={{\tau }_{1,klj}}-{{\gamma }^{{{h}_{k}}}}{{p}_{k,j+2}}{{\tau }_{1,kl,j+1}}={{\tau }_{1,klj}}-({{\tau }_{1,klj}}-{{\tau }_{1,kl,{{J}_{k}}-1}})={{\tau }_{1,kl,{{J}_{k}}-1}}.\]
Therefore, we only need to show that the following holds:
\begin{multline}
	\sum\limits_{kl}{{{Z}_{k}}{{z}_{kl}}}+Qq+\sum\limits_{tkl;2\le t\le {{T}_{k}}}{(\gamma {{\tau }_{t-1,kl,0}}-{{\rho }_{kl}}){{n}_{tkl}}}- \\ 
	\sum\limits_{kl}{{{\tau }_{1,kl,{{J}_{k}}-1}}\sum\nolimits_{j}{{{x}_{1,klj}}}}-\sum\limits_{kl}{({{\rho }_{kl}}-\gamma {{p}_{k,2}}{{\tau }_{{{h}_{k}},kl,1}}){{n}_{1,kl}}}\ge 0,\forall s\in S,\forall a\in {{A}_{s}}. \label{step1-main}
\end{multline}
Constraint \eqref{step1-main} holds if the minimum value of its left-hand side (LHS) is greater than or equal to zero. In this minimization problem, the variables are ${z_{kl}},\forall tkl$, $q$, ${{n}_{tkl}},\forall tkl$, and ${{x}_{1,klj}},\forall klj$, and they are subject to the state and action sets constrains, i.e. \eqref{eq:state}-\eqref{u2}. In the following, we show that the minimum of the LHS is always equal to zero. 
To minimize the LHS, we should have: 
\begin{enumerate}[label=(\alph*)]
	\item	${{z}_{kl}}=0,\forall kl$: ${{\zeta }^{z}}$ is much higher than ${{\zeta }^{q}}$, so it is not worth having a non-zero ${{z}_{kl}}$ to have a smaller $q$, or larger ${{x}_{1,klj}},\forall klj$ or ${{n}_{1,kl}},\forall tkl$.
	\item	$q\in \{0,2{{d}_{0,l}}\},\forall l$: If at least one visit has been assigned to day 1, all of them should be in one location $l$ for which we have $q=2{{d}_{0,l}}$. It is not worth having visits in a location ${l}'$ with ${{d}_{0,{l}'}}\le {{d}_{0,l}}$. Because they may result in $q>2{{d}_{0,l}}$, and also we have ${{\tau }_{1,k{l}',{{J}_{k}}-1}}\le {{\tau }_{1,kl,{{J}_{k}}-1}},\forall k$ based on (38).
	\item	$(\sum\nolimits_{kj}{{{x}_{1,klj}}}+\sum\nolimits_{k}{{{n}_{1,kl}}})\in \{0,\left\lfloor (\chi -2{{d}_{0,l}})/e \right\rfloor \}$: If the nurse is supposed to serve at least one visit in location $l$ on day 1, i.e. $q>0$, then capacity should be used to its maximum.
	\item	${{n}_{tkl}}=0,\forall tkl,2\le t\le {{T}_{k}}-1$: We know ${{\rho }_{kl}}=\gamma {{\tau }_{{{T}_{k}}-1,kl,0}},{{T}_{k}}\ge 2$, and therefore $\gamma {{\tau }_{t-1,kl,0}}-{{\rho }_{kl}}>0,\forall tkl,2\le t\le {{T}_{k}}-1$.  
	\item	${{n}_{1,kl}}=0,\forall kl$, $T_k \ge 2$: We have ${{\tau }_{1,kl,{{J}_{k}}-1}}\ge {{\rho }_{kl}}-\gamma p_{k,2} {{\tau }_{{{h}_{k}},kl,1}}$.
\end{enumerate}

Based on (d) and (e), we have ${{n}_{{{T}_{k}},kl}}={{y}_{kl}},\forall kl$. The minimum of the LHS obtained from (a)-(e) is equal to zero.
\subsection{Step 2}
Based on the explanations in Step 1 above, the hypothesized solution generates binding constraints in the ALP model only for state-action pairs for which ${{z}_{kl}}=0,\forall kl$, $\sum\limits_{kl}{{{\tau }_{1,kl,{{J}_{k}}-1}}\sum\nolimits_{j}{{{x}_{1,klj}}}}=Qq$, and ${{y}_{kl}}={{n}_{{{T}_{k}},kl}},\forall kl$. For sake of simplicity, we shrink the set of state-action pairs satisfying these conditions to the following one:
\begin{multline}
	{{\Lambda }^{*}}=\{(s,a)|{{z}_{kl}}=0,\forall kl, 
	{{x}_{1,klj}}\in \{0,\left\lfloor (\chi -2{{d}_{0,l}})/e \right\rfloor \},\forall klj, 
	{{x}_{tklj}}\in \{0,x^{max}_{tklj}\},\forall klj;t\ge 2, \\ 
	{{n}_{tkl}}=0,\forall tkl,t\le {{T}_{k}}-1,T_k \ge 2,  
	{{y}_{kl}}={{n}_{{{T}_{k}},kl}},\forall kl, n_{T_k,kl}\in \{0,y^{max}_{kl}\},\forall kl,T_k \ge 2, \\	n_{T_k,kl}\in \{0,\left\lfloor (\chi -2{{d}_{0,l}})/e \right\rfloor \},\forall kl\},T_k = 1,q \le \chi\}.  
\end{multline}
To prove the optimality of the hypothesized solution using complementary slackness, we must show the existence of a dual feasible solution $\beta^*$ for which $\beta^* (s,a)=0,\forall (s,a)\notin {{\Lambda }^{*}},s\in S,a\in {{A}_{s}}$, and all the constraints in the D-ALP model associated with non-zero ALP parameter values are tight. Let $\Lambda =\{(s,a)|\beta^* (s,a)>0\}$. Hence, $\Lambda \subset {{\Lambda }^{*}}$. Therefore, we need to find the necessary and sufficient conditions under which the following set of equations are satisfied:
\begin{subequations}
	\begin{align}
		& (1-\gamma )\sum\limits_{(s,a)\in \Lambda }{\beta^* (s,a)}=1, \label{CS_1} \\
		& \sum\limits_{(s,a)\in \Lambda }{\beta^* (s,a)\left( {{x}_{tklj}}-\gamma {{p}_{k,j+1}}({{x}_{1,kl,j-1}+{n}_{1,kl}}) \right)}=\EX_{\alpha}[{{X}_{tklj}}],\forall tklj,t={{h}_{k}},j = 1, \label{CS_2_2} \\
		& \sum\limits_{(s,a)\in \Lambda }{\beta^* (s,a)\left( {{x}_{tklj}}-\gamma {{p}_{k,j+1}}{{x}_{1,kl,j-1}} \right)}=\EX_{\alpha}[{{X}_{tklj}}],\forall tklj,t={{h}_{k}},j\ge 2, \label{CS_2} \\
		& \sum\limits_{(s,a)\in \Lambda }{\beta^* (s,a)\left( {{x}_{tklj}}-\gamma {{x}_{t+1,klj}} \right)}=\EX_{\alpha}[{{X}_{tklj}}],\forall tklj,t<{{h}_{k}},j\ge 1, \label{CS_3} \\
		& \sum\limits_{(s,a)\in \Lambda }{\beta^* (s,a)\left( {{x}_{tklj}}-\gamma \left( {{x}_{t+1,klj}}+{{n}_{t+1,kl}} \right) \right)}=\EX_{\alpha}[{{X}_{tklj}}],\forall tklj,t\le {{T}_{k}}-1,j=0, \label{CS_4} \\
		& \sum\limits_{(s,a)\in \Lambda }{\beta^* (s,a)({{y}_{kl}}-\gamma \EX[{{Y}_{kl}}])}=\EX_{\alpha}[{{Y}_{kl}}],\forall kl, \label{CS_5}
	\end{align}
\end{subequations}
where \eqref{CS_2_2} to \eqref{CS_5} are the binding versions of  \eqref{D_ALP_2} and \eqref{D_ALP_3}.
For sake of simplicity, let $x^{max}_{tklj}=x^{max},\forall tklj$, $y^{max}_{kl}=y^{max},\forall kl$, $\EX_{\alpha}[{{X}_{tklj}}]=\varepsilon,\forall tklj$, $\EX_{\alpha}[{{Y}_{kl}}]={{\lambda }_{kl}},\forall kl$, $x_{l}^{m}=\left\lfloor (\chi -2{{d}_{0,l}})/e \right\rfloor $, $\beta _{tklj}^{x}=\sum\limits_{(s,a)\in \Lambda ,{{x}_{tklj}}>0}{\beta^* (s,a)}$, $\beta _{tkl}^{n}=\sum\limits_{(s,a)\in \Lambda ,{{n}_{tkl}}>0}{\beta^* (s,a)}$ and $\beta _{kl}^{y}=\sum\limits_{(s,a)\in \Lambda ,{{y}_{kl}}>0}{\beta^* (s,a)}$.
\subsubsection{Equation \eqref{CS_5}}
We can write equation  \eqref{CS_5} as $\sum\limits_{(s,a)\in \Lambda }{\beta^* (s,a){{y}_{kl}}}=\frac{{{\lambda }_{kl}}}{1-\gamma },\forall kl$. Thus, we have:
\begin{equation}
	\beta _{kl}^{y}=\frac{{{\lambda }_{kl}}}{(1-\gamma )x_{l}^{m}},\forall kl,T_k=1, \label{Beta_y_1}
\end{equation}
\begin{equation}
	\beta _{kl}^{y}=\frac{{{\lambda }_{kl}}}{(1-\gamma )y^{max}},\forall kl,T_k \ge 2. \label{Beta_y_2}
\end{equation}
Consequently, we have:
\begin{equation}
	\beta _{1,kl}^{n}=\frac{{{\lambda }_{kl}}}{(1-\gamma )x_{l}^{m}},\forall kl,T_k = 1, \label{Beta_n_1}
\end{equation}
\begin{equation}
	\beta _{1,kl}^{n}=0,\forall kl,T_k \ge 2. \label{Beta_n_2}
\end{equation}

\subsubsection{Equation \eqref{CS_4}}
We can write equation \eqref{CS_4} as follows:
\begin{multline}
	\sum\limits_{(s,a)\in \Lambda }{\beta^* (s,a){{x}_{tklj}}}-\gamma \sum\limits_{(s,a)\in \Lambda }{\beta^* (s,a){{x}_{t+1,klj}}}-\gamma \sum\limits_{(s,a)\in \Lambda }{\beta^* (s,a){{n}_{t+1,kl}}}=\varepsilon, \\ 
	\forall tklj,t\le {{T}_{k}}-1,j=0. \label{CS_4_2}
\end{multline}
Then, we have:
\begin{equation}
	\beta _{tklj}^{x}=\frac{1}{{x^{max}}}\left( \sum\nolimits_{i=0}^{{{T}_{k}}-t-1}{{{\gamma }^{i}}} \right)\varepsilon+\frac{1}{{x^{max}}}\left( \frac{{{\gamma }^{{{T}_{k}}-t}}{{\lambda }_{kl}}}{1-\gamma } \right),\forall tkl,2\le t\le {{T}_{k}}-1,j=0. \label{CS_4_beta_1}
\end{equation}
\begin{equation}
	\beta _{tklj}^{x}=\frac{1}{x_{l}^{m}}\varepsilon+\frac{\gamma }{x_{l}^{m}}\left( \left( \sum\nolimits_{i=0}^{{{T}_{k}}-3}{{{\gamma }^{i}}} \right)\varepsilon+\frac{{{\gamma }^{{{T}_{k}}-2}}{{\lambda }_{kl}}}{1-\gamma } \right),\forall tklj,t=1,{{T}_{k}}\ge 2,j=0. \label{CS_4_beta_2_3}
\end{equation}

\subsubsection{Equation \eqref{CS_2_2}}
We can write equation \eqref{CS_2_2} as follows:
\[\sum\limits_{(s,a)\in \Lambda }{\beta^* (s,a){{x}_{tklj}}}-\gamma {{p}_{k,j+1}}\sum\limits_{(s,a)\in \Lambda }{\beta^* (s,a)({{x}_{1,kl,j-1}}+n_{1,kl})}=\varepsilon,\forall tklj,t={{h}_{k}},j= 1.\]
Therefore, we have:
\begin{equation}
	\beta _{tklj}^{x}=\frac{1}{x_{l}^{m}}\varepsilon+\gamma {{p}_{k,2}}(\beta _{1,kl,j-1}^{x}+\beta _{1,kl}^{n}),\forall tklj,t={{h}_{k}}=1,j=1, \label{CS_2_2_beta_1}
\end{equation}
\begin{equation}
	\beta _{tklj}^{x}=\frac{1}{{x^{max}}}\left( \varepsilon+x_{l}^{m}\gamma {{p}_{k,j+1}}(\beta _{1,kl,j-1}^{x}+\beta _{1,kl}^{n}) \right),\forall tklj,t={{h}_{k}},{{h}_{k}}\ge 2,j=1. \label{CS_2_2_beta_2}
\end{equation}

\subsubsection{Equation \eqref{CS_2}}
We can write equation \eqref{CS_2} as follows:
\[\sum\limits_{(s,a)\in \Lambda }{\beta^* (s,a){{x}_{tklj}}}-\gamma {{p}_{k,j+1}}\sum\limits_{(s,a)\in \Lambda }{\beta^* (s,a){{x}_{1,kl,j-1}}}=\varepsilon,\forall tklj,t={{h}_{k}},j\ge 2.\]
Therefore, we have:
\begin{equation}
	\beta _{tklj}^{x}=\frac{1}{x_{l}^{m}}\left( \sum\nolimits_{i=0}^{j-1}{{{\gamma }^{i}}\prod\limits_{{i}'=0}^{i-1}{{{p}_{k,j+1-{i}'}}}} \right)\varepsilon+{{\gamma }^{j}}\prod\limits_{{i}'=2}^{j+1}{{{p}_{k,{i}'}}}(\beta _{1,kl,0}^{x}+\beta _{1,kl}^{n}),\forall tklj,t={{h}_{k}}=1,j\ge 2. \label{CS_2_beta_1}
\end{equation}
\begin{equation}
	\beta _{tklj}^{x}=\frac{1}{{x^{max}}}\left( \varepsilon+x_{l}^{m}\gamma {{p}_{k,j+1}}\beta _{1,kl,j-1}^{x} \right),\forall tklj,t={{h}_{k}},{{h}_{k}}\ge 2,j\ge 2. \label{CS_2_beta_2}
\end{equation}

\subsubsection{Equation \eqref{CS_3}}
We can write equation \eqref{CS_3} as follows:
\[\sum\limits_{(s,a)\in \Lambda }{\beta^* (s,a){{x}_{tklj}}}-\gamma \sum\limits_{(s,a)\in \Lambda }{\beta^* (s,a){{x}_{t+1,klj}}}=\varepsilon,\forall tklj,t<{{h}_{k}},j\ge 1.\]
Thus, we have:
\begin{equation}
	\beta _{tklj}^{x}=\frac{1}{x_{l}^{m}}\left( \varepsilon+\gamma x^{max}\beta _{t+1,klj}^{x} \right),\forall tklj,t=1,{{h}_{k}}=2,j\ge 1. \label{CS_3_beta_1}
\end{equation}
\begin{equation}
	\beta _{tklj}^{x}=\frac{1}{{x^{max}}}\left( \varepsilon\sum\nolimits_{i=0}^{{{h}_{k}}-t-1}{{{\gamma }^{i}}}+{{\gamma }^{{{h}_{k}}-t}}x^{max}\beta _{{{h}_{k}},klj}^{x} \right),\forall tklj,2\le t<{{h}_{k}},{{h}_{k}}\ge 3,j\ge 1. \label{CS_3_beta_2}
\end{equation}
\begin{equation}
	\beta _{tklj}^{x}=\frac{1}{x_{l}^{m}}\left( \varepsilon+\gamma x^{max}\beta _{t+1,klj}^{x} \right),\forall tklj,t=1,{{h}_{k}}\ge 3,j\ge 1. \label{CS_3_beta_3}
\end{equation}

In summary, $\beta _{kl}^{y},\forall kl$, $\beta _{1,kl}^{n},\forall kl$ and $\beta _{tklj}^{x},\forall tklj$  can be calculated using equations \eqref{Beta_y_1}-\eqref{Beta_n_2}, \eqref{CS_4_beta_1}, and \eqref{CS_4_beta_2_3}-\eqref{CS_3_beta_3}. Because these formulas are dependent to each other, they should be calculated in the following order. First, $\beta _{tklj}^{x},\forall tkl,j=0$ using \eqref{CS_4_beta_1} and \eqref{CS_4_beta_2_3}. Second, $\beta _{tklj}^{x},\forall tklj,t={{h}_{k}}=1,j\ge 1$ using \eqref{CS_2_2_beta_1} and \eqref{CS_2_beta_1}. Third, $\beta _{tklj}^{x},\forall kl,t={{h}_{k}}=2$ using \eqref{CS_2_2_beta_2} and \eqref{CS_2_beta_2} and $\beta _{tklj}^{x},\forall kl,t=1,{{h}_{k}}=2$ using \eqref{CS_3_beta_1} sequentially for $j$ equal to 1 to ${{J}_{k}}-1$. Finally, $\beta _{tklj}^{x},\forall kl,t={{h}_{k}}$ using \eqref{CS_2_beta_2}, $\beta _{tklj}^{x},\forall kl,2\le t\le {{h}_{k}}-1$ using \eqref{CS_3_beta_2}, and $\beta _{tklj}^{x},\forall kl,t=1$ using \eqref{CS_3_beta_3} sequentially for $j$ equal to 1 to ${{J}_{k}}-1$.
\subsection{Step 3}
Now, we only need to show that there exists a dual solution satisfying the necessary and sufficient condition in Step 2 above. It can be easily shown that such a solution $\hat \beta^*$ can be described as follows: 
\begin{enumerate}[label=(\alph*)]
	\item For $\forall klj$, consider a $(s,a)\in \Lambda $ in which the only non-zero state variable is ${{x}_{1,klj}}=x_{l}^{m}$, and let $\hat \beta^* (s,a)=\beta _{1,klj}^{x}$.
	\item	For $\forall klj,T_k=1$, consider a $(s,a)\in \Lambda $ in which the only non-zero state variable is ${{n}_{1,kl}}=x_{l}^{m}$, and let $\hat \beta^* (s,a)=\beta _{1,kl}^{n}$.
	\item	For $\forall tklj,t\ge 2$, consider a $(s,a)\in \Lambda $ in which the only non-zero state variable is ${{x}_{tklj}}=x^{max}$($x^{max}$ is very large), and let $\hat \beta^* (s,a)=\beta _{tklj}^{x}$. 
	\item	For $\forall kl,T_k \ge 2$, consider a $(s,a)\in \Lambda $ in which the only non-zero state variable is ${{y}_{kl}}=y^{max}$($y^{max}$ is very large), and let $\hat \beta^* (s,a)=\beta _{kl}^{y}$.
	\item	Let $\hat \beta^* (s=0,a=0)={{\beta }^{0}}=\frac{1}{1-\gamma }-\sum\nolimits_{tklj}{\beta _{tklj}^{x}}+\sum\nolimits_{kl,T_k = 1}{\beta _{1,kl}^{n}}+\sum\nolimits_{kl,T_k \ge 2}{\beta _{kl}^{y}}$.
	\item For all other $(s,a)\in \Lambda $, let $\hat \beta^* (s,a)=0$.
\end{enumerate}

In conclusion, the hypothesized primal solution and a corresponding dual solution would satisfy the complementary slackness conditions if having the $\beta _{tklj}^{x},\forall tklj$, $\beta _{1,kl}^{n},\forall kl$ and $\beta _{kl}^{y},\forall kl$ values obtained from \eqref{CS_2_2}-\eqref{CS_5}, equation \eqref{CS_1} holds, i.e.,
$\sum\limits_{(s,a)\in \Lambda }{\hat \beta^* (s,a)}=\frac{1}{(1-\gamma )}$. It can be easily shown that $\sum\limits_{(s,a)\in \Lambda }{\hat \beta^* (s,a)}= \sum\nolimits_{tklj}{\beta _{tklj}^{x}}+\sum\nolimits_{kl,T_k = 1}{\beta _{1,kl}^{n}}+\sum\nolimits_{kl,T_k \ge 2}{\beta _{kl}^{y}}+{\beta }^{0}$, where $\underset{x^{max}\to \infty }{\mathop{\lim }}\,\beta _{tklj}^{x}=0,\forall tklj,t\ge 2$ and $\underset{y^{max}\to \infty }{\mathop{\lim }}\,\beta _{kl}^{y}=0,\forall kl,T_k \ge 2$. Because $\sum\nolimits_{klj}{\beta _{1,klj}^{x} +\sum\nolimits_{kl}{\beta _{1,kl}^{n} }<\frac{1}{1-\gamma }}$ is satisfied (as a condition of the theorem), $\sum\limits_{(s,a)\in \Lambda }{\hat \beta^* (s,a)}=\frac{1}{(1-\gamma )}$ holds with a non-negative $\beta^0$.

\section{Analysis of the heuristic reduction techniques for the ALP}\label{sec:effects_of_accel}
Assume a $2\times3$ rectangular area consisting of 6 regions and one service type. We determined the ALP parameters for 25 random instances of the same problem setting using each of the ALP, ALP-2I, ALP-1D, and ALP-1D-2I with average execution times of 3751.72, 11.44, 13.93, and 0.66. The random instances were generated by changing the service type characteristics, $\bar{D}$, overtime limit, and $\varepsilon $. The application of the heuristic reduction techniques does not result in any loss of optimality and decreases the execution time of obtaining the ALP parameters significantly for the instances considered. The best version is the ALP-1D-2I model taking advantage of both techniques mentioned in Section \ref{sec:Acceleration-techniques-for}.

\section{Additional information about the numerical experiments}\label{sec:E_add_num}
In Table \ref{K=2-instances}, we provide the input parameters of the random instances with two service types whose results were reported in Table \ref{table-K=2}.

\begin{table}[h]
	\centering
	\caption{Characteristics of the service types for instances with $K=2$}
	\footnotesize
	\fontsize{5}{7}\selectfont
	\begin{tabular}{M{1.5cm}M{1.5cm}M{1.5cm}M{1.5cm}M{1.5cm}M{1.5cm}}
		\hline
		Ins \# & $k$ & ${{h}_{k}}$ & $\EX [\hat{J}_k]$ & ${{e}_{k}}$ & ${{T}_{k}}$ \\ \hline
		\multirow{2}{*}{1} &  1 & 7 & 10 & 0.50 & 1                \\ 
		&     	2 & 2 & 6 & 0.50 & 3             \\ \hline
		\multirow{2}{*}{2} &  1 & 3 & 6 & 0.25 & 4               \\ 
		&     	2 & 4 & 10 & 0.50 & 2              \\ \hline
		\multirow{2}{*}{3} &  1 & 7 & 6 & 0.75 & 2                \\ 
		&     	2 & 3 & 6 & 1.00 & 4              \\ \hline
		\multirow{2}{*}{4} &  1 & 2 & 6 & 1.00 & 2               \\ 
		&     	2 & 1 & 6 & 0.25 & 1              \\ \hline
		\multirow{2}{*}{5} &  1 & 2 & 10 & 0.25 & 1               \\ 
		&     	2 & 7 & 8 & 0.75 & 1              \\ \hline
		\multirow{2}{*}{6} &  1 & 2 & 10 & 1.00 & 4               \\ 
		&     	2 & 1 & 6 & 0.75 & 1              \\ \hline
		\multirow{2}{*}{7} &  1 & 5 & 6 & 0.50 & 3                \\ 
		&     	2 & 6 & 6 & 1.00 & 4              \\ \hline
		\multirow{2}{*}{8} &  1 & 5 & 8 & 0.75 & 1                \\ 
		&     	2 & 1 & 8 & 0.25 & 4              \\ \hline
		\multirow{2}{*}{9} &  1 & 5 & 10 & 0.75 & 5               \\ 
		&     	2 & 7 & 6 & 1.00 & 4              \\ \hline
		\multirow{2}{*}{10} &  1 & 5 & 10 & 1.00 & 2                \\ 
		&     	2 & 7 & 6 & 0.25 & 1              \\ \hline
		\multirow{2}{*}{11} & 1 & 3 & 8 & 0.75 & 5                  \\ 
		& 2 & 2 & 8 & 0.25 & 3                    \\ \hline
		\multirow{2}{*}{12} &   1 & 4 & 6 & 0.25 & 4		                 \\ 
		&   2 & 2 & 8 & 1.00 & 1              \\ \hline
		\multirow{2}{*}{13} &  1 & 7 & 8 & 1.00 & 5            \\ 
		&    		2 & 7 & 6 & 1.00 & 2                   \\ \hline
		\multirow{2}{*}{14} &    1 & 6 & 10 & 1.00 & 1              \\ 
		&      2 & 4 & 6 & 0.25 & 2             \\ \hline
		\multirow{2}{*}{15} &  1 & 3 & 10 & 1.00 & 2                \\ 
		&      2 & 4 & 6 & 0.75 & 3             \\ \hline
		\multirow{2}{*}{16} &    1 & 1 & 8 & 1.00 & 1             \\ 
		&    2 & 1 & 6 & 0.75 & 3                \\ \hline
		\multirow{2}{*}{17} &     1 & 6 & 8 & 0.25 & 2              \\ 
		&       	2 & 1 & 10 & 0.50 & 1           \\ \hline
		\multirow{2}{*}{18} &  1 & 7 & 8 & 0.50 & 1                \\ 
		&    2 & 6 & 10 & 0.25 & 4               \\ \hline
		\multirow{2}{*}{19} &    1 & 3 & 8 & 0.25 & 4               \\ 
		&       2 & 5 & 8 & 0.50 & 4           \\ \hline
		\multirow{2}{*}{20} &  1 & 3 & 10 & 0.25 & 2                \\ 
		&     	2 & 1 & 10 & 0.50 & 1              \\ \hline
		\multirow{2}{*}{21} &  1 & 6 & 6 & 0.75 & 5               \\ 
		&     	2 & 2 & 6 & 0.75 & 1              \\ \hline
		\multirow{2}{*}{22} &  1 & 7 & 10 & 0.50 & 4                \\ 
		&     	2 & 4 & 10 & 1.00 & 2              \\ \hline
		\multirow{2}{*}{23} &  1 & 5 & 6 & 1.00 & 1                \\ 
		&     	2 & 4 & 6 & 0.25 & 4              \\ \hline
		\multirow{2}{*}{24} &  1 & 1 & 10 & 0.25 & 4                \\ 
		&     	2 &4 & 8 & 0.75 & 4              \\ \hline
		\multirow{2}{*}{25} &  1 & 6 & 6 & 1.00 & 1               \\ 
		&     	2 & 1 & 10 & 0.25 & 3              \\ \hline
	\end{tabular}\label{K=2-instances}
\end{table}

Table \ref{table-simulation-multiple} presents the results for a single random instance for each $K\in \{3,4,6,8,10\}$ with default values for the input parameters (reported in Table \ref{table-general-settings} for $K>1$). The ALP policy outperforms the SB policy in the instances with 3, 4, and 8 service types, and has a bit better average percentage gap (that is not statistically significant) in the other two.
\begin{table}[h]    
	\caption{Simulation results for instances with multiple service types, $K\in \{3,4,6,8,10\}$, and overtime}
	\centering
	\footnotesize
	\resizebox{\textwidth}{!}{%
		\fontsize{9}{13}\selectfont
		\begin{tabular}{M{1cm}M{1cm}M{1cm}M{1cm}M{1cm}M{1cm}M{1cm}M{1cm}M{1cm}M{1cm}M{1cm}M{1cm}M{1cm}M{1cm}M{1cm}M{1cm}M{2cm}M{2cm}}
			\hline
			\multirow{2}{*}{$K$} &
			\multicolumn{3}{c}{Rejection hours} &
			\multicolumn{3}{c}{Diversion hours} &
			\multicolumn{3}{c}{Overtime (hours)} &
			\multicolumn{3}{c}{Travel time (hours)} &
			\multicolumn{3}{c}{Tour length (hours)} &
			\multicolumn{2}{c}{Gap\% (SD)} \\ \cline{2-18}
			& ALP  & Myopic & SB & ALP  & Myopic & SB & ALP  & Myopic& SB & ALP  & Myopic& SB & ALP  & Myopic& SB & ALP & SB              \\ \hline
			3  & 1.25 & 0.11 & 1.04 & 0.31 & 1.26 & 0.59 & 1.46 & 1.60 & 1.47 & 1.46 & 1.56 & 1.67 & 9.39 & 9.59 & 9.45 & \textbf{21.70 (1.73)} & 10.62 (1.49)
			\\  
			4  & 1.15 & 0.08 & 0.96 & 0.34 & 1.25 & 0.60 & 1.53 & 1.69 & 1.55 & 1.54 & 1.63 & 1.78 & 9.48 & 9.69 & 9.54 & \textbf{21.60 (1.69)} & 12.07 (1.48)
			\\  
			6  & 1.06 & 0.33 & 1.49 & 0.52 & 1.17 & 0.34 & 1.47 & 1.56 & 1.29 & 1.49 & 1.58 & 1.72 & 9.37 & 9.54 & 9.26 & \textbf{15.31 (1.93)} & \textbf{13.77 (2.05)}
			\\  
			8  & 1.25 & 0.35 & 1.09 & 0.36 & 1.08 & 0.54 & 1.46 & 1.63 & 1.47 & 1.53 & 1.61 & 1.75 & 9.38 & 9.62 & 9.46 & \textbf{14.57 (1.20)} & 10.43 (1.33)
			\\  
			10 & 1.48 & 0.30 & 1.87 & 0.34 & 1.31 & 0.25 & 1.40 & 1.55 & 1.15 & 1.51 & 1.58 & 1.78 & 9.28 & 9.53 & 9.10 & \textbf{16.98 (1.17)} & \textbf{15.47 (1.27)}
			\\ \hline
		\end{tabular}\label{table-simulation-multiple}%
	}
\end{table}

\section{A heuristic for the multi-nurse problem}\label{sec:E_multi-nurse}
An important direction for future research would be to extend the proposed ALP approach to solve the multi-nurse version of the problem, which involves multiple nurses with different depots. There are multiple service types, and each nurse may be eligible to serve a subset of them. The objective would be to satisfy the demand for home care services as much as possible while minimizing the costs. It is possible not to employ some of the available nurses. A simple heuristic is stated in the following to solve this problem. In this heuristic, the multi-nurse problem is decomposed to multiple single-nurse problems while they affect each other. The ALP parameters of the corresponding problem to each nurse are obtained using the following algorithm: \vspace{5pt}

\noindent
\textbf{Algorithm EC.2. Obtaining the ALP parameters for each nurse}
\begin{steps}
	\item Create a sequence $\Psi=(\Psi_1,...,\Psi_N)$ of the available $N$ nurses, and let $i=1$. 
	\item Get the ALP parameters for $\Psi_i$ using the ALP-1D-2I. Calculate the average rejection rate for each service type and radius (i.e., regions with the same distance from $\Psi_i$'s depot) using simulation. If it is worth having $\Psi_{i}$, go to the next step; otherwise, go to Step 4.
	\item Update the referral rate of each service type and region for the remaining nurses by subtracting the amount covered by $\Psi_i$ from the current value. 
	\item If $i<N$, let $i=i+1$ and back to Step 2; otherwise, stop.
\end{steps}

The sequence of the available nurses in Step 1 may affect the quality of the final solution. Different strategies can be used to determine this sequence. For example, we can cluster the demand in the whole area and assign priorities to nurses based on their closeness to the densest clusters. Then, we start from a high-priority nurse close to the densest cluster, and proceed to the next high-priority nurse that is closest to the previously investigated one. An example of a cluster-based heuristic for the VRP problem can be found in \citep{Dondo2007}.

In the single-nurse problem we solve in Step 2, we only consider the service types for which $\Psi_{i}$ is eligible. Moreover, only regions around the $\Psi_{i}$'s depot that are accessible for this nurse are taken into account. The accessibility of regions for a specific nurse can be defined in different ways. We can define a region as accessible for a nurse if she/he can perform at least one visit of any service type, she/he is eligible for, considering the travel times and limited shift length. However, if this definition includes many regions and makes the problem intractable, we can consider some other definitions such as only regions located inside a specific area diameter around the depot (like one hour) are accessible. In Step 2, the nurse $\Psi_{i}$ is employed if the value of covering additional demand by her/him is more than the cost of employing her/him.

In Step 3, the way we update the remaining referral rates may affect the final solution. It can be done symmetrically or asymmetrically. For example, suppose that the average rejection rate of a specific radius is 50\% for a specific nurse. In a symmetric update, the remaining referral rates in all regions located at that radius may become a half of their current values. However, in an asymmetric one, those of one side of the nurse's depot may become equal to zero and those of its opposite side that are closer to another nurse coming after her/him in $\Psi$ remain unchanged.

Once the ALP parameters corresponding to each nurse are available, an action set can be generated each day for each of them one by one according to their order in $\Psi$.

\end{document}